\documentclass[11pt, reqno]{amsart}

\usepackage{amssymb}

\evensidemargin\oddsidemargin

\begin{document}
\baselineskip=18pt
\setcounter{page}{1}

\newcommand{\eqnsection}{
\renewcommand{\theequation}{\thesection.\arabic{equation}}
    \makeatletter
    \csname  @addtoreset\endcsname{equation}{section}
    \makeatother}
\eqnsection
   

\def\a{\alpha}
\def\B{{\bf B}}
\def\cA{{\mathcal{A}}} 
\def\cD{{\mathcal{D}}} 
\def\cG{{\mathcal{G}}} 
\def\cH{{\mathcal{H}}} 
\def\cL{{\mathcal{L}}} 
\def\cI{{\mathcal{I}}}
\def\CC{{\mathbb{C}}} 
\def\Dap{{\rm D}_{0+}^\a} 
\def\Dm{{\rm D}_{-}^\a} 
\def\Dp{{\rm D}_{+}^\a} 
\def\Ea{E_\a}
\def\esp{{\mathbb{E}}} 
\def\F{{\bf F}}
\def\Farl{{\F}_{\a,\lbd,\rho}}
\def\bF{{\bar F}}
\def\bG{{\bar G}}
\def\G{{\bf G}}
\def\Ga{{\Gamma}} 
\def\Gal{{\G_{\a,\lbd}}} 
\def\GG{{\bf \Gamma}}
\def\ii{{\rm i}} 
\def\Iap{{\rm I}_{0+}^\a} 
\def\Im{{\rm I}_{-}^\a} 
\def\Ip{{\rm I}_{+}^\a} 
\def\L{{\bf L}}
\def\lbd{\lambda}
\def\lacc{\left\{}
\def\lcr{\left[}
\def\lpa{\left(}
\def\lva{\left|}
\def\M{{\bf M}}
\def\NN{{\mathbb{N}}} 
\def\pb{{\mathbb{P}}}
\def\QQ{{\mathbb{Q}}} 
\def\R{{\bf R}}
\def\rl{{\mathbb{R}}}
\def\racc{\right\}}
\def\rpa{\right)}
\def\rcr{\right]}
\def\rva{\right|}
\def\sga{\sigma^{(\a)}}
\def\T{{\bf T}}
\def\Un{{\bf 1}}
\def\X{{\bf X}}
\def\Y{{\bf Y}}
\def\W{{\bf W}}
\def\Z{{\bf Z}}
\def\Warl{{\W}_{\a,\lbd,\rho}}
\def\Za{{\Z_\a}}

\def\elaw{\stackrel{d}{=}}
\def\claw{\stackrel{d}{\longrightarrow}}
\def\elaw{\stackrel{d}{=}}
\def\qed{\hfill$\square$}


\newtheorem{prop}{Proposition}[section]
\newtheorem{coro}[prop]{Corollary}
\newtheorem{lem}[prop]{Lemma}
\newtheorem{ex}[prop]{Example}
\newtheorem{exs}[prop]{Examples}
\newtheorem{rem}[prop]{Remark}
\newtheorem{theo}[prop]{Theorem}
\newtheorem{defi}[prop]{Definition}
\newtheorem{conj}[prop]{Conjecture}

\newcommand{\noi}{\noindent}
\newcommand{\dis}{\displaystyle }

\title{Some properties of the Kilbas-Saigo function}

\author[L.~Boudabsa]{Lotfi Boudabsa}

\address{Institut de Math\'ematiques, Ecole Polytechnique F\'ed\'erale de Lausanne, CH-1015 Lausanne, Switzerland. {\em Email}: {\tt lotfi.boudabsa@epfl.ch}}

\author[T.~Simon]{Thomas Simon}

\address{Laboratoire Paul Painlev\'e, UMR 8524, Universit\'e de Lille,  Cit\'e Scientifique, F-59655 Villeneuve d'Ascq Cedex, France. {\em Email}: {\tt thomas.simon@univ-lille.fr}}

\keywords{Complete monotonicity; Convex ordering; Double Gamma function; Kilbas-Saigo function; Le Roy function; Mittag-Leffler function; Stable subordinator}

\subjclass[2010]{33B15; 33E12; 60E15; 60G52}

\begin{abstract} 
 
We characterize the complete monotonicity of the Kilbas-Saigo function on the negative half-line. We also provide the exact asymptotics at $-\infty$, and uniform hyperbolic bounds are derived. The same questions are addressed for the classical Le Roy function. The main ingredient for the proof is a probabilistic representation of these functions in terms of the stable subordinator.

\end{abstract}

\maketitle

\section{Introduction}

\label{FurthA}

The Kilbas-Saigo function is a three-parameter entire function having the convergent series representation
$$E_{\a,m,l} (z)\; =\; 1\; +\; \sum_{n\ge 1} \lpa\prod_{k=1}^n \frac{\Ga(1+ \a((k-1)m +l))}{\Ga(1 + \a((k-1)m + l +1))}\rpa z^n, \qquad z\in\CC,$$
where the parameters are such that for $\a,m >0$ and $l > -1/\a.$ It can be viewed as a generalization of the one- or two-parameter Mittag-Leffler function since, with the standard notation for the latter,
$$E_{\a,1,0}(z)\; =\; \sum_{n\ge 0} \frac{z^n}{\Ga(1 +\a n)}\; =\; \Ea(z)$$
and
$$E_{\a,1,\frac{\beta-1}{\a}}(z)\; =\; \Ga(\beta) \sum_{n\ge 0} \frac{z^n}{\Ga(\beta +\a n)}\; =\;\Ga(\beta) E_{\a,\beta} (z)$$ 
for every $\a,\beta > 0$ and $z\in\CC.$ This function was introduced in \cite{KS95} as the solution to some integro-differential equation with Abelian kernel on the half-line, and we refer to Chapter 5.2 in \cite{GKMR} for a more recent account, including an extension to complex values of the parameter $l.$ In our previous paper \cite{BSV}, written in collaboration with P.~Vallois, it was shown that certain Kilbas-Saigo functions are moment generating functions of Riemannian integrals of the stable subordinator. This observation made it possible to define rigorously some Weibull and Fr\'echet distributions of fractional type via an independent exponential random variable and the stable subordinator - see \cite{BSV} for details. In the present paper, we wish to take the other way round and use the probabilistic connection in order to deduce some non-trivial analytical properties of the Kilbas-Saigo function. 

In Section 2, we tackle the problem of the complete monotonicity on the negative half-line. This problem dates back to Pollard in 1948 for the one-parameter Mittag-Leffler function - see e.g. Section 3.7.2 in \cite{GKMR} for details and references. It was shown in \cite{BSV} that  for every $m >0$ and $\a\in (0,1]$ the function $x\mapsto E_{\a,m,m-1}(-x)$ is completely monotone, extending Pollard's result and solving an open problem stated in \cite{DMV}. In Theorem \ref{KSCM} below, we characterize the complete monotonicity of $x\mapsto E_{\a,m,l}(-x)$ by $\a\in [0,1]$ and $l \ge m-1/\a.$ We also give an explicit representation, albeit complicated in general, of the underlying positive random variable. Along the way, we study an interesting family of Mellin transforms given as the quotient of four double Gamma functions.

In Section 3, we establish uniform hyperbolic bounds on the negative half-line for two families of completely monotonic Kilbas-Saigo functions, extending the bounds obtained in \cite{TS2014} for the classical Mittag-Leffler function. The argument in \cite{TS2014} relied on stochastic and convex orderings and was rather lenghty. We use here the same kind of arguments, but the proof is shorter and more transparent thanks to the connection with the stable subordinator. The latter also implies some monotonicity properties on $m\mapsto E_{\a,m,l} (x)$ for every $x\in\rl$ - see Proposition \ref{Mono} below.

In Section 4, we address the question of the asymptotic behaviour at $-\infty$ in the completely monotonic case $\a\in (0,1]$ and $l \ge m-1/\a.$ It is shown in Theorem 5.5 of \cite{GKMR} that in the general case $\a, m > 0$ and $l > m-1/\a,$ the entire function $E_{\a,m,l}(z)$ has order $\rho = 1/\a$ and type $\sigma = 1/m.$ However, precise asymptotics along given directions of the complex plane do not seem to have been investigated as yet, as is the case - see e.g. Proposition 3.6 in \cite{GKMR} - for the classical Mittag-Leffler function. For the negative half-line and $\a\in (0,1]$, the asymptotics are different depending on whether $l=m+1/\a$ or $l > m+ 1/\a.$ In the former case, the behaviour is in $c_{\a, m}\, x^{-(1+1/m)}$ with a non-trivial constant $c_{\a, m}$ obtained from the connection with the fractional Fr\'echet distribution and given in terms of the double Gamma function - see Proposition \ref{MASF} and Remark \ref{rF2} (c) below. In the latter case, the behaviour is in $c_{\a, m, l}\, x^{-1}$ with a uniform speed and a simple constant $c_{\a, m, l}$ given in terms of the standard Gamma function - see Proposition \ref{ASKS} below. The method for the case $l > m+ 1/\a$ relies on the computation of the Mellin transform of the positive function $E_{\a,m,l}(-x),$ which is obtained from the proof of its complete monotonicity, and is interesting in its own right - see Remark \ref{RCM} (c) below. Along the way, we provide the exact asymptotics of the fractional Weibull and Fr\'echet densities at both ends of their support and we give a series of probabilistic factorizations. The latter enhance the position of the fractional Fr\'echet distribution, which is in one-to-one correspondence with the boundary Kilbas-Saigo function $E_{\a, m, m-1/\a} (x),$ as an irreducible factor - see Remark \ref{rF2} (a) below.    
 
In the last Section 5, we pay attention to the so-called Le Roy function with parameter $\a > 0.$ The latter is a simple generalization of the exponential function defined by
$$\cL_\a(z)\; =\; \sum_{n\ge 0} \frac{z^n}{(n!)^\a}, \qquad z\in\CC.$$
Introduced in \cite{L1899} in the context of analytic continuation, a couple of years before the Mittag-Leffler function, the Le Roy function has been much less studied. It was shown in \cite{BSV} that this function encodes for $\a\in[0,1]$ a Gumbel distribution of fractional type, as the moment generating function of the perpetuity of the $\a-$stable subordinator. This fact is recalled in Proposition \ref{LRCM} below, together with a characterization of the moment generating property. The exact asymptotic behaviour at $-\infty$ is also derived for $\a\in [0,1],$ completing the original result of Le Roy. Finally, the non-increasing character of $\a\mapsto \cL_\a(x)$ on $[0,1]$ for every $x\in\rl$ is established by convex ordering. It is worth mentioning that the latter property is an open problem - see Conjecture \ref{MLMono} below - for the Mittag-Leffler function.

As in \cite{BSV}, an important role is played throughout the paper by Barnes' double Gamma function $G(z;\delta)$ which is the unique solution to the functional equation $G(z+1;\delta) = \Gamma(z\delta^{-1}) G(z;\delta)$ with normalization $G(1;\delta) =1,$ and its associated Pochhammer type symbol
$$[a;\delta]_s\; =\; \frac{G(a +s;\delta)}{G(a;\delta)}\cdot$$
We have gathered in an Appendix all the needed facts and formul\ae\, on this double Gamma function, whose connection with the Kilbas-Saigo function has probably a broader focus than the content of the present paper (we leave this topic open to further research). 

\section{Complete monotonicity on the negative half-line} In this section, we wish to characterize the property that the function $x\mapsto E_{\a,m,l}(-x)$  is completely monotone (CM) on $(0,\infty).$ We begin with the following result on the above generalized Pochhammer symbols, which is reminiscent of Proposition 5.1 and Theorem 6.2 in \cite{D2010} and has an independent interest.

\begin{lem}
\label{Gtype}
Let $a,b,c,d$ and $\delta$ be positive parameters. There exists a positive random variable $Z =\Z[a,c;b,d;\delta]$ such that 
\begin{equation}
\label{Gdelta}
\esp[Z^s]\; =\; \frac{[a;\delta]_s[c;\delta]_s}{[b;\delta]_s[d;\delta]_s}
\end{equation}
for every $s > 0,$ if and only if $b+d \le a+c$ and $\inf\{b,d\} \le \inf\{a, c\}.$ This random variable is absolutely continuous on $(0,\infty)$, except in the degenerate case $a=b=c=d.$ Its support is $[0,1]$ if $b+d = a+c$ and $[0,\infty)$ if $b+d < a+c.$ 
\end{lem}

\proof

We giscard the degenerate case $a=b=c=d,$ which is obvious with $Z=1.$ By (\ref{Malm}) and some rearrangements - see also (2.15) in \cite{LS}, we first rewrite
$$\log\lpa \frac{[a;\delta]_s[c;\delta]_s}{[b;\delta]_s[d;\delta]_s}\rpa \; =\; \kappa\, s\; + \; \int_{-\infty}^0 (e^{sx} - 1 - sx) \lpa \frac{e^{- b\vert x\vert} + e^{- d\vert x\vert} - e^{- a\vert x\vert} - e^{- c\vert x\vert}}{\vert x\vert (1-e^{-\vert x\vert})(1-e^{-\delta \vert x\vert})}\rpa dx$$
for every $s > 0,$ where $\kappa$ is some real constant. By convexity, it is easy to see that if $b+d \le a+c$ and $\inf\{b,d\} \le \inf\{a, c\},$ then the function $z\mapsto z^b + z^d - z^a - z^c$ is positive on $(0,1).$ This implies that the function
$$x\;\mapsto\;  \frac{e^{- b\vert x\vert} + e^{- d\vert x\vert} - e^{- a\vert x\vert} - e^{- c\vert x\vert}}{\vert x\vert (1-e^{-\vert x\vert})(1-e^{-\delta \vert x\vert})}$$
is positive on $(-\infty, 0)$ and that it can be viewed as the density of some L\'evy measure on $(-\infty,0),$ since it integrates $1\wedge x^2.$ By the L\'evy-Khintchine formula, there exists a real infinitely divisible random variable $Y$ such that 
$$\esp[e^{sY}] \; =\;\frac{[a;\delta]_s[c;\delta]_s}{[b;\delta]_s[d;\delta]_s}$$  
for every $s > 0,$ and the positive random variable $Z = e^Y$ satisfies (\ref{Gdelta}). Since we have excluded the degenerate case, the L\'evy measure of $Y$ is clearly infinite and it follows from Theorem 27.7 in \cite{S1999} that $Y$ has a density and the same is true for $Z.$

Assuming first $b+d = a+c,$ a Taylor expansion at zero shows that the density of the L\'evy measure of $Y$ integrates $1\wedge \vert x\vert$ and we deduce from (\ref{Malm}) the simpler formula
$$\log\esp[e^{sY}]\; =\; \log\lpa \frac{[a;\delta]_s[c;\delta]_s}{[b;\delta]_s[d;\delta]_s}\rpa \; =\; -\int^{\infty}_0 (1 - e^{-sx}) \lpa \frac{e^{- bx} + e^{- dx} - e^{- ax} - e^{- cx}}{x(1-e^{-x})(1-e^{-\delta x})}\rpa dx.$$
By the L\'evy-Khintchine formula, this shows that the ID random variable $Y$ is negative. Moreover, its support is $(-\infty,0]$ since its L\'evy measure has full support and its drift coefficient is zero - see Theorem 24.10 (iii) in \cite{S1999}, so that the support of $Z$ is $[0,1].$ 

Assuming second $b+d < a+c,$ the same Taylor expansion as above shows that the density of the L\'evy measure of $Y$ does not integrate $1\wedge \vert x\vert$ and the real L\'evy process associated to $Y$ is hence of type C with the terminology of \cite{S1999} - see Definition 11.9 therein. By Theorem 24.10 (i) in \cite{S1999}, this implies that $Y$ has full support on $\rl,$ and so does $Z$ on $\rl^+\!.$ \\
  
It remains to prove the only if part of the Lemma. Assuming $a \le d$ and $b\le c$ without loss of generality, we first observe that if $a < b$ then the function
$$s\;\mapsto\; \frac{[a;\delta]_s[c;\delta]_s}{[b;\delta]_s[d;\delta]_s}$$
is real-analytic on $(-b,\infty)$ and vanishes at $s = - a > -b,$ an impossible property for the Mellin transform of a positive random variable. The necessity of $b+d\le a+c$ is slightly more subtle and hinges again upon infinite divisibility. First, setting $\varphi(z) = z^b + z^d - z^a - z^c$ and $z_* =\inf\{ z > 0, \; \varphi(z) < 0\},$ it is easy to see by convexity and a Taylor expansion at 1 that if $b+d > a+c,$ then $z_* < 1$ and $\varphi(z) < 0$ on $(z_*,1)$ with $\varphi(z)\sim (b+d-a-c)(z-1)$ as $z\to 1.$ Introducing next the ID random variable $V$ with Laplace exponent
$$\log\esp[e^{sV}]\; =\; - \kappa\, s\; + \; \int_{\log z_*}^0 (e^{sx} - 1 - sx) \lpa \frac{e^{- a\vert x\vert} + e^{- c\vert x\vert} - e^{- b\vert x\vert} - e^{- d\vert x\vert}}{\vert x\vert (1-e^{-\vert x\vert})(1-e^{-\delta \vert x\vert})}\rpa dx,$$
we obtain the decomposition
$$\log\lpa \frac{[a;\delta]_s[c;\delta]_s}{[b;\delta]_s[d;\delta]_s}\rpa\; +\; \log\esp[e^{sV}] \; =\;\int_{-\infty}^{\log z_*} (e^{sx} - 1 - sx) \lpa \frac{e^{- b\vert x\vert} + e^{- d\vert x\vert} - e^{- a\vert x\vert} - e^{- c\vert x\vert}}{\vert x\vert (1-e^{-\vert x\vert})(1-e^{-\delta \vert x\vert})}\rpa dx,$$
whose right-hand side is the Laplace exponent of some ID random variable $U$ having an atom because its L\'evy measure, whose support is bounded away from zero, is finite - see Theorem 27.4 in \cite{S1999}. On the other hand, the random variable $V$ has an absolutely continous and infinite L\'evy measure and hence it has also a density. If there existed $Z$ such that (\ref{Gdelta}) holds, then the independent decomposition $U\elaw V  + \log Z$ would imply by convolution that $U$ has a density as well. This contradiction finishes the proof of the Lemma.

\endproof

\begin{rem} {\em (a) By the Mellin inversion formula, the density of $\Z[a,c;b,d;\delta]$ is expressed as
$$f(x) \; =\; \frac{1}{2\ii\pi x} \int_{s_0 -\ii \infty}^{s_0 + \ii\infty} x^{-s} \lpa \frac{[a;\delta]_s[c;\delta]_s}{[b;\delta]_s[d;\delta]_s} \rpa \, ds$$
over $(0,\infty)$ for any $s_0 > -\inf\{b,d\}.$ From this expression, it is possible to prove that this density is real-analytic over the interior of the support. We omit details. Let us also mention by Remark 28.8 in \cite{S1999} that this density is positive over the interior of its support.

\medskip

(b) With the standard notation for the Pochhammer symbol, the aforementioned Proposition 5.1 and Theorem 6.2 in \cite{D2010} show that
$$s\;\mapsto\;\frac{(a)_s(c)_s}{(b)_s(d)_s}$$
is the Mellin transform of a positive random variable if and only if $b+d \ge a+c$ and $\inf\{b,d\} \ge \inf\{a, c\}.$ This fact can be proved exactly as above, in writing
$$\log\lpa \frac{(a)_s(c)_s}{(b)_s(d)_s}\rpa \; =\; - \int^{\infty}_0 (1-e^{-sx}) \lpa \frac{e^{- ax} + e^{- c x} - e^{- b x} - e^{- d x}}{x(1-e^{-x})}\rpa dx.$$
This expression also shows that the underlying random variable has support $[0,1]$ and that it is absolutely continuous, save for $a+c = b+d$ where it has an atom at zero. We refer to \cite{D2010} for an exact expression of the density on $(0,1)$ in terms of the classical hypergeometric function.}

\end{rem}

We can now characterize the CM property for $E_{\a,m,l}(-x)$ on $(0,\infty).$

\begin{theo}
\label{KSCM}
Let $\a, m > 0$ and $l > -1/\a.$ The Kilbas-Saigo function
$$x\;\mapsto\; E_{\a,m,l}(-x)$$
is {\em CM} on $(0,\infty)$ if and only if $\a\le 1$ and $l\ge m-1/\a.$ Its Bernstein representation is  
\begin{equation}
\label{BKS}
E_{\alpha,m,l}(-x)\; =\; \esp\lcr \exp -x \lacc \X_{\a,m,l}\;\times\int_0^\infty \lpa 1+ \sga_t\rpa^{-\a (m+1)} \, dt\racc\rcr
\end{equation}
with $\delta =1/\a m$ and $\X_{\a,m,l} = \Z[1+1/m, (\a l +1)\delta;1, 1/m + (\a l +1)\delta;\delta].$ 
\end{theo}

\proof
Assume first $\a \le 1$ and $l\ge m-1/\a$ and let 
$$\Y_{\a,m,l}\;=\; \X_{\a,m,l}\;\times\int_0^\infty \lpa 1+ \sga_t\rpa^{-\a (m+1)}\, dt.$$ 
By Proposition 2.4 in \cite{LS}, and Lemma \ref{Gtype}, its Mellin transform is
\begin{eqnarray*}
\esp[(\Y_{\a,m,l})^s] & = &  \delta^{s}\, \frac{[1+\delta;\delta]_s[(\a l +1)\delta;\delta]_s}{[1;\delta]_s[1/m +(\a l +1)\delta;\delta]_s}\\
& = & \Ga(1+s)\,\times\frac{[(\a l +1)\delta;\delta]_s}{[1/m +(\a l+1)\delta;\delta]_s}
\end{eqnarray*}
where in the second equality we have used (\ref{Conca4}). By Fubini's theorem, the moment generating function of $\Y_{\a,m,l}$ reads
\begin{eqnarray*}
\esp[e^{z \Y_{\a,m,l}}] & = & \sum_{n\ge 0}\; \esp[(\Y_{\delta,m,\eta})^n]\,\frac{z^n}{n!}\\
& = & \sum_{n\ge 0} \lpa\frac{[(\a l +1)\delta;\delta]_n}{[1/m +(\a l +1)\delta;\delta]_n}\rpa z^n\\
& = &  \sum_{n\ge 0} \lpa\prod_{j=0}^{n-1}\frac{\Ga(\a(jm +l)+1)}{\Ga(\a(jm +l +1)+1)}\rpa z^n\; =\; E_{\a,m,l}(z)
\end{eqnarray*}
for every $z\ge 0,$ where in the third equality we have used (\ref{Conca1}) repeatedly. The latter identity is extended analytically to the whole complex plane and we get, in particular,
$$ E_{\a,m,l}(-x) \; =\; \esp[e^{-x \Y_{\a,m,l}}],\qquad x\ge 0.$$
This shows that $E_{\a,m,l}(-x)$ is CM with the required Bernstein representation.\\

We now prove the only if part. If $E_{\a,m,l}(-x)$ is CM, then we see by analytic continuation that $E_{\a,m,l}(z)$ is the moment generating function on $\CC$ of the underlying random variable $X,$ whose positive integer moments read 
$$\esp[X^n]\; =\; n!\,\times\lpa\prod_{j=0}^{n-1}\frac{\Ga(\a(jm +l)+1)}{\Ga(\a(jm +l +1)+1)}\rpa,\quad n\ge 0.$$ 
If $\a > 1,$ Stirling's formula implies $\esp[X^n]^{\frac{1}{n}}\to 0$ as $n\to\infty$ so that $X\equiv 0,$ a contradiction because $E_{\a,m,l}$ is not a constant. If $\a = 1$ and $l+1 < m,$ then
$$\esp[X^n]\; = \;\frac{n!}{(c)_n m^n}\; \sim\; \frac{n^{1-c}}{m^n}\quad \mbox{as $n\to\infty,$}$$ 
with $c = (l+1)/m \in (0,1).$ In particular, the Mellin transform $s\mapsto\esp[X^s]$ is analytic on $\{\Re(s)\ge 0\},$ bounded on $\{\Re(s) = 0\},$ and has at most exponential growth on $\{\Re(s) > 0\}$ because
$$\vert \esp[X^s]\vert\; \le\; \esp\lcr X^{\Re(s)}\rcr\; =\; \lpa \esp\lcr X^{[\Re(s)] +1}\rcr\rpa^{\frac{\Re(s)}{[\Re(s)] +1}}$$
by H\"older's inequality. On the other hand, the Stirling type formula (\ref{Stirling}) implies, after some simplifications, 
$$\delta^{-s} \frac{[1+\delta;\delta]_s[c;\delta]_s}{[1;\delta]_s[c+\delta;\delta]_s}\; = \; \delta^{-s} s^{1-c} (1+ o(1))\quad \mbox{as $\vert s\vert \to\infty$ with $\vert \arg s\vert < \pi$}$$
and this shows that the function on the left-hand side, which is analytic on $\{\Re(s)\ge 0\},$ has at most linear growth on $\{\Re(s) = 0\}$ and at most exponential growth on $\{\Re(s) > 0\}.$ Moreover, the above analysis clearly shows that
$$\esp[X^n] \; = \; \delta^{-n} \frac{[1+\delta;\delta]_n[c;\delta]_n}{[1;\delta]_n[c+\delta;\delta]_n}$$
for all $n\ge 0$ and by Carlson's theorem - see e.g. Section 5.81 in \cite{T1939}, we must have
$$\esp[X^s] \; = \; \delta^{-s} \frac{[1+\delta;\delta]_s[c;\delta]_s}{[1;\delta]_s[c+\delta;\delta]_s}$$
for every $s > 0,$ a contradiction since Lemma \ref{Gtype} shows that the right-hand side cannot be the Mellin transform of a positive random variable if $c < 1.$ The case $\a < 1$ and  $l+1/\a < m$ is analogous. It consists in identifying the bounded sequence
$$ \frac{1}{n!}\,\times\lpa\prod_{j=0}^{n-1}\frac{\Ga(\a(jm +l +1)+1)}{\Ga(\a(jm +l)+1)}\rpa$$
as the values at non-negative integer points of the function
$$\delta^{-s}\times \frac{[1;\delta]_s[1/m + (\a l +1)\delta;\delta]_s}{[1+\delta;\delta]_s[(\a l +1)\delta;\delta]_s}\; =\; \delta^{-s} e^{-(1-\a) s \ln (s) + \kappa s + O(1)}\quad\mbox{as $\vert s\vert \to \infty$ with $\vert \arg s\vert < \pi,$}$$
where the purposeless constant $\kappa$ can be evaluated from (\ref{Stirling}). On $\{\Re(s) \ge 0\},$ we see that this function has growth at most $e^{\pi (1-\a) \vert s\vert /2}$ and we can again apply Carlson's theorem. We leave the details to the interested reader.

\endproof

\begin{rem}

\label{RCM}

{\em (a) When $m =1,$ applying (\ref{Conca1}) we see that the random variable $\X_{\a, 1,l}$ has Mellin transform
$$\esp[(\X_{\a,1,l})^s]\; =\; \frac{[2;\delta]_s[l+1/\a;\delta]_s}{[1;\delta]_s[1 + l + 1/\a;\delta]_s}\; =\; \frac{(\a)_{\a s}}{(\beta)_{\a s}}$$
with $\beta = 1+\a l\ge \a.$ This shows $\X_{\a, 1,l}\elaw \B_{\a,\beta -\a}^\a$ where $\B_{a,b}$ denotes, here and throughout, a standard Beta random variable with parameters $a,b >0.$ We hence recover the Bernstein representation of the CM function $\Ga(\beta) E_{\a,\beta}(-x)$ which was discussed in Remark 3.3 (c) in \cite{BSV}. Notice also the very simple expression for the Mellin transform
$$\esp[(\Y_{\a,1,l})^s] \; =\; \frac{\Ga(1+\a l)\Ga(1+s)}{\Ga(1 +\a (l+s))}\cdot$$

(b) Another simplification occurs when $l+1/\a = k m$ for some integer $k\ge 1.$ One finds
$$\esp[(\X_{\a,m,km -1/\a})^s]\; =\; \frac{[k;\delta]_s[1 +1/m;\delta]_s}{[1;\delta]_s[k +1/m;\delta]_s}\; =\; \prod_{j=1}^{k-1} \lpa \frac{(\a j m)_{u}}{(\a(jm+1))_{u}}\rpa$$
for $u =\a m s \ge 0,$ which implies
$$\X_{\a, m, km - 1/\a}\; \elaw\; \lpa \B_{\a m, \a}\times\cdots\times \B_{\a m (k-1), \a}\rpa^{\a m}.$$
In general, the law of the absolutely continuous random variable $\X_{\a,m,l}$ valued in $[0,1]$ seems to have a complicated expression.

\medskip

(c) As seen during the proof, the random variable $\Y_{\a,m,l}$ defined by the Bersntein representation 
$$E_{\a, m, l}(-x) \; =\; \esp [e^{-x \Y_{\a, m, l}}]$$
has Mellin transform 
\begin{equation}
\label{MellY}
\esp[(\Y_{\a,m,l})^s]\; =\; \Ga(1+s)\, \times\frac{[(\a l +1)\delta;\delta]_s}{[1/m +(\a l+1)\delta;\delta]_s}
\end{equation}
with $\delta = 1/\a m,$ for every $s > -1.$ By Fubini's theorem, this implies the following exact computation, which seems unnoticed in the literature on the Kilbas-Saigo function.
\begin{equation}
\label{Marr}
\int_0^\infty E_{\a, m, l}(-x)\, x^{s-1}\, dx \; = \; \Ga(s)\,\esp [\Y_{\a,m,l}^{-s}]\; = \; \Ga(s)\Ga(1-s)\, \times\frac{[(\a l +1)\delta;\delta]_{-s}}{[1/m +(\a l+1)\delta;\delta]_{-s}}
\end{equation}
for every $s\in (0,1).$ For $m = 1,$ we recover from (\ref{Conca1}) the formula
$$\int_0^\infty E_{\a, \beta}(-x)\, x^{s-1}\, dx \; = \; \frac{1}{\Ga(\beta)} \int_0^\infty E_{\a,1, \frac{\beta-1}{\a}}(-x)\, x^{s-1}\, dx \; = \;\frac{\Ga(s)\Ga(1-s)}{\Ga(\beta - \a s)}$$ 
which is given in (4.10.3) of \cite{GKMR}, as a consequence of the Mellin-Barnes representation of $E_{\a,\beta}(z).$ Notice that there is no such Mellin-Barnes representation for $E_{\a,m,l}(z)$ in general.}

\end{rem}

\section{Uniform hyperbolic bounds}

In Theorem 4 of \cite{TS2014}, the following uniform hyperbolic bounds are obtained for the classical Mittag-Leffler function:
\begin{equation}
\label{MLB}
\frac{1}{1 +\Ga(1-\a) x}\; \le\; E_\a(-x)\; \le\; \frac{1}{1 +\frac{1}{\Ga(1+\a)}\, x}
\end{equation}
for every $\a\in[0,1]$ and $x\ge 0.$ The constants in these inequalities are optimal because of the asymptotic behaviours
$$E_\a(-x)\; \sim\;\frac{1}{\Ga(1-\a) x}\quad\mbox{as $x\to\infty$}\qquad\mbox{and}\qquad 1- E_\a(-x)\; \sim\; \frac{x}{\Ga(1+\a)} \quad\mbox{as $x\to 0.$}$$ See \cite{M} and the references therein for some motivations on these hyperbolic bounds. In this section, we shall obtain analogous bounds for $E_{\a,m,m-1}(-x)$ and $E_{\a,m,m-\frac{1}{\a}}(-x)$ with $\a\in[0,1], m > 0.$ Those peculiar functions are associated to the fractional Weibull and Fr\'echet distributions defined in \cite{BSV}. Specifically, we will use the following representations as a moment generating function, obtained respectively in (3.1) and (3.4) therein:
\begin{equation}
\label{KS1}
E_{\alpha,m,m-1}(z)\; =\; \esp\lcr\exp \lacc z\, \int_0^\infty \lpa 1- \sga_t\rpa_+^{\a(m-1)} \, dt\racc\rcr
\end{equation}
and 
\begin{equation}
\label{KS2}
E_{\alpha,m,m-\frac{1}{\a}}(z)\; =\; \esp\lcr\exp \lacc z\, \int_0^\infty \lpa 1 + \sga_t\rpa^{-\a (m+1)} \, dt\racc\rcr
\end{equation}
for every $z\in\CC,$ where $\{\sga_t\, t\ge 0\}$ is the $\a-$stable subordinator normalized such that
$$\esp[e^{-\lbd \sga_t}]\; =\; e^{-t\lbd^\a},\qquad \lbd, t\ge 0.$$ Observe that these two formul\ae\, specify the general Bernstein representation (\ref{BKS}) in terms of the $\a-$stable subordinator only. We begin with the following monotonicity properties, of independent interest.
\begin{prop}
\label{Mono}
Fix $\a\in (0,1]$ and $x\in\rl.$ The functions 
$$m\; \mapsto\; E_{\a,m,m-1}(x)\qquad\mbox{and}\qquad m\; \mapsto\; E_{\a,m,m-\frac{1}{\a}}(x)$$
are decreasing on $(0,\infty)$ if $x > 0$ and increasing on $(0,\infty)$ if $x < 0.$  
\end{prop}

\proof

This is a consequence of (\ref{KS1}) resp. (\ref{KS2}), and the fact that $\sga_t > 0$ for every $t > 0.$

\endproof

\begin{rem}
{\em It would be interesting to know if the same property holds for $m \mapsto E_{\a,m,m - l}(x)$ and any $l \le 1/\a.$ In the case $l \not\in\{1, 1/\a\},$ this would require from (\ref{BKS}) a monotonicity analysis of the mapping $m\mapsto \X_{\a, m, m-l},$ which does not seem easy at first sight.}
\end{rem} 

As in \cite{TS2014}, our analysis to obtain the uniform bounds will use some notions of stochastic ordering. Recall that if $X,Y$ are real random variables such that $\esp[\varphi(X)] \le \esp[\varphi(Y)]$ for every $\varphi :\rl \to\rl$ convex, then $Y$ is said to dominate $X$ for the convex order, a property which we denote by $X\prec_{cx} Y.$ Another ingredient in the proof is the following infinite independent product

$$\T(a,b,c)\; =\; \prod_{n\ge 0} \lpa\frac{a+nb +c}{a+nb}\rpa \B_{a+ nb, c}.$$
We refer to Section 2.1 in \cite{LS} for more details on this infinite product, including the fact that it is a.s. convergent for every $a,b,c > 0.$ We also mention from Proposition 2 in \cite{LS} that its Mellin transform is 
$$\esp[\T(a,b,c)^s]\;=\; \lpa \frac{\Ga(ab^{-1})}{\Ga ((a+c)b^{-1})}\rpa^s\times\; \frac{[a+c; b]_s}{[a; b]_s}$$
for every $s > -a.$ The following simple result on convex orderings for the above infinite independent products has an independent interest.

\begin{lem}
\label{Tcx}
For every $a,b,c > 0$ and $d \ge c,$ one has
$$\T(a,b,c)\;\prec_{cx}\; \T(a,b,d).$$ 
\end{lem}

\proof

By the definition of $\T(a,b,c)$ and the stability of the convex order by mixtures - see Corollary 3.A.22 in \cite{SSh}, it is enough to show
$$(a+b) \B_{a,b}\;\prec_{cx}\;(a+c) \B_{a,c}$$ 
for every $a,b > 0$ and $c \ge b.$ Using again Corollary 3.A.22 in \cite{SSh} and the standard identity $\B_{a,c}\,\elaw\,\B_{a,b}\times\B_{a+b,c-b},$ we are reduced to
$$\lpa\frac{a+b}{a+c}\rpa\; =\; \esp[\B_{a+b,c-b}]\;\prec_{cx}\;\B_{a+b,c-b}$$  
which is a consequence of Jensen's inequality.

\endproof

The following result is a generalization of the inequalities (\ref{MLB}), which deal with the case $m=1$ only, to all Kilbas-Saigo functions $E_{\a,m,m-1}(-x).$ The argument is considerably simpler than in the original proof of (\ref{MLB}).

\begin{theo}
\label{KSO1}
For every $\a\in[0,1], m >0$ and $x\ge 0,$ one has
$$\frac{1}{1 +\Ga(1-\a) x}\; \le\; E_{\a,m,m-1}(-x)\; \le\; \frac{1}{1 +\frac{\Ga(1+\a(m-1))}{\Ga(1+\a m)}\, x}\cdot$$
\end{theo}

\proof

The first inequality is a consequence of Proposition \ref{Mono}, which implies in letting $m\to 0$ 
\begin{eqnarray*}
E_{\a,m,m-1}(-x) & \ge & \esp\lcr\exp \lacc -x\, \int_0^\infty \lpa 1- \sga_t\rpa_+^{-\a} \, dt\racc\rcr\\
& = & \esp\lcr e^{-x\,\Ga(1-\a)\,\L}\rcr\; = \; \frac{1}{1 +\Ga(1-\a) x}
\end{eqnarray*}
for $x\ge 0,$ where the first equality follows from Theorem 1.2 (b) (ii) in \cite{LS}. For the second inequality, we come back to the infinite product representation
\begin{eqnarray*}
\int_0^\infty \lpa 1- \sga_t\rpa_+^{\rho -\a} \, dt & \elaw & \frac{\Ga(\rho +1-\a)}{\Ga(\rho +1)}\; \T(1,\rho^{-1}, (1-\a)\rho^{-1})
\end{eqnarray*}
which follows from Theorem 1.2 (b) (i) in \cite{LS}, exactly as in the proof of Theorem 1.1 in \cite{BSV}. Lemma \ref{Tcx} implies then 
\begin{eqnarray*}
\int_0^\infty \lpa 1- \sga_t\rpa_+^{\rho -\a} \, dt & \prec_{cx} & \frac{\Ga(\rho +1-\a)}{\Ga(\rho +1)}\; \T(1,\rho^{-1}, \rho^{-1})\; \elaw\; \frac{\Ga(\rho +1-\a)}{\Ga(\rho +1)}\;\L
\end{eqnarray*}
where the identity in law follows from (2.7) in \cite{LS}. Using (\ref{KS1}) with $\rho = \a m$ and the convexity of $t\mapsto e^{-xt}$, we obtain the required
$$E_{\a,m,m-1}(-x)\; \le\; \frac{1}{1 +\frac{\Ga(1+\a(m-1))}{\Ga(1+\a m)}\, x}\cdot$$
\endproof

\begin{rem}
\label{rKS1}
{\em (a) As for the classical case $m=1$, these bounds are optimal because of the asymptotic behaviours
$$1\,-\, E_{\a, m, m-1}(-x)\; \sim\; \frac{\Ga(1+\a (m-1))}{\Ga(1+\a m)}\, x \quad\mbox{as $x\to 0$}$$
and
$$E_{\a, m, m-1}(-x)\; \sim\;\frac{1}{\Ga(1-\a) x}\quad\mbox{as $x\to\infty.$}$$
The behaviour at zero is plain from the definition, whereas the behaviour at infinity will be given after Remark \ref{rW2} below.

\medskip

(b) It is easy to check that the above proof also yields the upper bound
$$E_{\a,m,m-1}(x)\; \le\; \frac{1}{(1 -\Ga(1-\a) x)_+}$$
for every $\a\in[0,1], m > 0$ and $x\ge 0,$ which seems unnoticed even in the classical case $m=1.$
}
\end{rem}

Our next result is a uniform hyperbolic upper bound for the Kilbas-Saigo function $ E_{\a, m, m-\frac{1}{\a}}(-x),$ with a power exponent which will be shown to be optimal in Remark \ref{rF2} (c) below, and also an optimal constant because 
$$1\,-\, E_{\a, m, m-\frac{1}{\a}}(-x)\; \sim\; \lpa 1 +\frac{1}{m}\rpa \times \frac{\Ga(1+\a m)\, x}{\Ga(1+\a(m+1))}\quad \mbox{as $x\to 0.$}$$ 
\begin{prop}
\label{KSO2}
For every $\a\in(0,1], m >0$ and $x\ge 0,$ one has
$$E_{\a,m,m-\frac{1}{\a}}(-x)\; \le\; \frac{1}{\lpa 1 + \frac{\Ga(1+\a m)}{\Ga(1+\a (m+1))}\, x\rpa^{1+\frac{1}{m}}}\cdot$$ 
\end{prop}

\proof

The inequality is derived by convex ordering as in Theorem \ref{KSO1}: setting, here and throughout, $\GG_a$ for a Gamma random variable with parameter $a >0,$ one has
\begin{eqnarray*}
\int_0^\infty \lpa 1+ \sga_t\rpa_+^{-\rho -\a} \, dt & \elaw & \frac{\Ga(\rho)}{\Ga(\rho +\a)}\; \T(1 +\a\rho^{-1},\rho^{-1}, (1-\a)\rho^{-1})\\
& \prec_{cx} & \frac{\Ga(\rho)}{\Ga(\rho +\a)}\; \T(1 +\a\rho^{-1},\rho^{-1}, \rho^{-1})\; \elaw\; \frac{\Ga(\rho +1)}{\Ga(\rho +1+\a)}\;\GG_{1+\frac{\a}{\rho}}
\end{eqnarray*}
where the first identity follows from Corollary 3 in \cite{LS} as in the proof of Theorem 1.1 in \cite{BSV}, the convex ordering from Lemma \ref{Tcx} and the second identity from (2.7) in \cite{LS}. Then, using (\ref{KS2}) with $\rho = \a m,$ we get the required inequality.

\endproof

As in Theorem \ref{KSO1}, we believe that there is also a uniform lower bound, with a more complicated optimal constant which can be read off from the asymptotic behaviour of the density at zero obtained in Proposition \ref{MASF} below:

\begin{conj}
\label{Cjlow}
For every $\a\in(0,1], m >0$ and $x\ge 0,$ one has
\begin{equation}
\label{ConjLow}
E_{\a,m,m-\frac{1}{\a}}(-x) \;\ge\; \frac{1}{(1 +(\a m)^{-\frac{\a}{m+1}}(\Ga(1+\a)\, G(1-\a; \a m)\, G(1+\a; \a m))^{-\frac{m}{m+1}}\, x)^{1+\frac{1}{m}}}\cdot
\end{equation}
\end{conj}
Unfortunately, the proof of this general inequality still eludes us. The monotonicity property observed in Proposition \ref{Mono} does not help here, giving only the trivial lower bound zero. The discrete factorizations which are used in \cite{TS2014} are also more difficult to handle in this context, because the Mellin transform underlying $E_{\a,m,m-\frac{1}{\a}}$ is expressed in terms of generalized Pochhammer symbols. In the case $m=1,$ we could however get a proof of (\ref{ConjLow}). The argument relies on the following representation, observed in Remarks 3.1 (d) and 3.3 (c) of \cite{BSV}:
\begin{equation}
\label{T1a}
E_{\alpha,1,1-\frac{1}{\a}}(z)\; =\;\Ga(\a)E'_{\a,\a}(z)\; =\;\Ga(1+\a)E'_{\a}(z)\; =\; \Ga(1+\a)\esp \lcr T_\a\, e^{z T_\a}\rcr\; =\;  \esp \lcr e^{z T_\a^{(1)}}\rcr
\end{equation}
for every $z\in\CC,$ where $T_\a = \inf\{ t > 0,\; \sga_t > 1\}$ is the first-passage time above one of the $\a-$stable subordinator and $T^{(1)}_\a$ its usual size-bias of order one. 

\begin{prop}
\label{KSO3}
For every $\a\in(0,1)$ and $x\ge 0,$ one has
$$E_{\a,1,1-\frac{1}{\a}}(-x)\; \ge\; \frac{1}{\lpa 1 +\sqrt{\frac{\Ga(1-\a)}{\Ga(1+\a)}}\, x\rpa^{\! 2}}\cdot$$
\end{prop}

\proof By (\ref{T1a}) and since
$$\esp \lcr e^{-x \,\GG_2}\rcr\; =\; \frac{1}{\lpa 1 + x\rpa^2}$$
for every $x \ge 0,$ it is enough to show, reasoning exactly as in the proof of Theorem 4 in \cite{TS2014}, that
\begin{equation}
\label{StoO}
T_\a^{(1)}\,\prec_{st}\, \sqrt{\frac{\Ga(1-\a)}{\Ga(1+\a)}}\, \GG_2,
\end{equation}
where $\prec_{st}$ stands for the usual stochastic order between two real random variables. Recall that $X\prec_{st} Y$ means $\pb[X\ge x]\le\pb[Y\ge x]$ for every $x\in\rl.$ Since $T_{1/2}\elaw 2\sqrt{\GG_{1/2}},$ the case $\a =1/2$ is explicit and the stochastic ordering can be obtained directly. More precisely, the densities of both random variables in (\ref{StoO}) are respectively given by
$$\frac{x}{2}\; e^{-x^2/4}\qquad\mbox{and}\qquad \frac{x}{2}\; e^{-x/\sqrt{2}}$$
on $(0,\infty),$ where they cross only once at $x= 2\sqrt{2}.$ It is a well-known and easy result that this single intersection property yields (\ref{StoO}) - see Theorem 1.A.12 in \cite{SSh}.

The argument for the case $\a\neq 1/2$ is somehow analogous, but the details are more elaborate because the density of $T_\a^{(1)}$ is not explicit anymore. We proceed as in Theorem C of \cite{TS2014} and first consider the case where $\a$ is rational. Setting $\a = p/n$ with $n > p$ positive integers and $X_\a =T_\a^{(1)}$ we have, on the one hand,
\begin{eqnarray*}
\esp[(X_{\a})^{n s}] & = & \frac{\esp[(T_\a)^{1+ns}]}{\esp[T_\a]} \\ 
& = & \frac{\Ga(2+ns) \Ga(1+pn^{-1})}{\Ga(1+pn^{-1} +ps)}\\
& = & \frac{n^{ns}}{p^{ps}}\,\times \,\esp\lcr \lpa\B_{\frac{2}{n}, \frac{1}{p} -\frac{1}{n}}\rpa^s\rcr\,\times\,\frac{\prod_{i=3}^{n+1} (in^{-1})_s}{\prod_{j=2}^{p} (jp^{-1} +n^{-1})_s}
\end{eqnarray*}
for every $s > -2n^{-1},$ where we have used the well-known identity $T_\a\elaw (\sga_1)^{-\a}$ in the second equality, whereas in the third equality we have used repeatedly the Legendre-Gauss multiplication formula for the Gamma function - see e.g. Theorem 1.5.2 in \cite{AAR1999}. The same formula implies, on the other hand,
\begin{eqnarray*}
\esp\lcr\lpa\sqrt{\frac{\Ga(1-\a)}{\Ga(1+\a)}}\; \GG_2\rpa^{\!\! ns}\, \rcr & = & \frac{n^{ns}\, \kappa_\a^s}{p^{ps}}\,\times \,\esp\lcr \lpa\GG_{\frac{2}{n}}\rpa^s\rcr\,\times\,\lpa \prod_{i=3}^{n+1} (in^{-1})_s\rpa\\
& = & \frac{n^{ns}}{p^{ps}}\,\times\,\esp\lcr \lpa\kappa_\a\,\times \,\GG_{\frac{2}{n}}\,\times\, \prod_{j=2}^{p}\; \GG_{\frac{j}{p} +\frac{1}{n}}\rpa^{\!s}\,\rcr\,\times\, \frac{\prod_{i=3}^{n+1} (in^{-1})_s}{\prod_{j=2}^{p} (jp^{-1} +n^{-1})_s}
\end{eqnarray*}
for every $s > -2n^{-1},$ with the notation
$$\kappa_\a\; =\; \lpa\prod_{i=1}^p\frac{\Ga(ip^{-1} - n^{-1})}{\Ga(ip^{-1} + n^{-1})}\rpa^{\frac{n}{2}}.$$ Since 
$$\frac{\prod_{i=3}^{n+1} (in^{-1})_s}{\prod_{j=2}^{p} (jp^{-1} +n^{-1})_s}\; =\;\esp\lcr \lpa\prod_{i=2}^{p}\; \B_{\frac{i+1}{n}, \frac{i}{p} - \frac{i}{n}}\,\times\! \prod_{j=p+1}^{n} \GG_{\frac{j+1}{n}}\rpa^{\!s}\;\rcr$$
for every $s > -3n^{-1},$ by factorization and Theorem 1.A.3(d) in \cite{TS2014} we are finally reduced to show
$$\B_{\frac{2}{n}, \frac{1}{p} -\frac{1}{n}} \prec_{st}  \lpa\prod_{i=1}^p\;\frac{\Ga(ip^{-1} - n^{-1})}{\Ga(ip^{-1} + n^{-1})}\rpa^{\frac{n}{2}} \!\times \,\GG_{\frac{2}{n}}\,\times\, \prod_{j=2}^{p}\; \GG_{\frac{j}{p} +\frac{1}{n}}$$
for every $n >p$ positive integers. The latter is equivalent to 
$$(\B_{\frac{2}{n}, \frac{1}{p} -\frac{1}{n}})^{\frac{2}{n}} \prec_{st}  \lpa\prod_{i=2}^p\;\frac{\Ga(ip^{-1} - n^{-1})}{\Ga(ip^{-1} + n^{-1})} \rpa\times \lpa\GG_{\frac{2}{n}}\,\times\, \prod_{j=2}^{p}\; \GG_{\frac{j}{p} +\frac{1}{n}}\rpa^{\frac{2}{n}}$$
and this is proved via the single intersection property exactly as for (5.1) in \cite{TS2014}: the random variable on the left-hand side has an increasing density on $(0,1),$ whereas the random variable on the right-hand side has a decreasing density on $(0,\infty),$ both densities having the same positive finite value at zero. We omit details. This completes the proof of (\ref{StoO}) when $\a$ is rational. The case when $\a$ is irrational follows then by a density argument.

\endproof

\begin{rem}
\label{rKS3}
{\em It is easy to check from (\ref{Conca2}) and (\ref{Rhorho}) that
$$\frac{\Ga(1+\a)}{\Ga(1-\a)}\; =\;\a^\a\,\Ga(1+\a)\, G(1-\a; \a )\, G(1+\a; \a),$$
so that Proposition \ref{KSO3} leads to (\ref{ConjLow}) for $m=1,$ in accordance with the estimate (\ref{Fala}). In general, the absence of a tractable complement formula for the product $G(1-\a; \delta )\, G(1+\a; \delta)$ makes however the constant in (\ref{ConjLow}) more difficult to handle.}
\end{rem}

Our last result in this section gives optimal uniform hyperbolic bounds for the generalized Mittag-Leffler functions $E_{\a,\beta}(-x)$ whenever they are completely monotone, that is for $\beta \ge \a$ - see the above Remark \ref{RCM} (a). This can be viewed as another generalization of (\ref{MLB}).

\begin{prop}
\label{KSO4}
For every $\a\in(0,1], \beta > \a$ and $x\ge 0,$ one has 
$$\frac{1}{\lpa 1 +\sqrt{\frac{\Ga(1-\a)}{\Ga(1+\a)}}\, x\rpa^2}\;\le\; \Ga(\a)\,E_{\a,\a}(-x)\; \le\; \frac{1}{\lpa 1 +\frac{\Ga(1+\a)}{\Ga(1+ 2\a)}\, x\rpa^2}$$
and
$$\frac{1}{1 +\frac{\Ga(\beta-\a)}{\Ga(\beta)}\, x}\;\le\; \Ga(\beta)\,E_{\a,\beta}(-x)\; \le\; \frac{1}{1 +\frac{\Ga(\beta)}{\Ga(\beta+\a)}\, x}\cdot$$
\end{prop}

\proof The bounds for $E_{\a,\a}(-x)$ are a direct consequence of (\ref{T1a}), Proposition \ref{KSO2} and Proposition \ref{KSO3}. Notice that letting $\a\to 1$ leads to the trivial bound $0\le e^{-x}\le (2/(2+x))^2.$ To handle the bounds for $\beta > \a,$ we first recall from Remark \ref{RCM} (a) that
$$\Ga(\beta)\,E_{\a,\beta}(-x)\; =\; \Ga(\beta)\,E_{\a,1,\frac{\beta -1}{\a}}(-x)\; =\; \esp\lcr e^{-x \, \Y_{\a,1,l}}\rcr$$
with $l = (\beta - 1)/\a > 1-1/\a$ and $\Y_{\a, 1, l} \elaw \B_{\a,\beta -\a}^\a \times\, T_\a^{(1)}.$ Moreover, one has 
\begin{equation}
\label{MomXab}
\esp \lcr (\Y_{\a,1,l})^s\rcr\; =\; \frac{\Ga(1+s) \Ga(\beta)}{\Ga(\beta + \a s)}
\end{equation}
for every $s > -1,$ which implies the factorization $\L \elaw \Y_{\a,1,l}\times (\GG_\beta)^{\a}.$ Since, by Jensen's inequality,
$$\frac{\Ga(\beta +\a)}{\Ga(\beta)}\; =\; \esp\lcr (\GG_\beta)^\a\rcr \prec_{cx} (\GG_\beta)^{\a},$$
we deduce from Corollary 3.A.22 in \cite{SSh} the convex ordering
$$\Y_{\a,1,l}\prec_{cx} \frac{\Ga(\beta)}{\Ga(\beta +\a)}\;\L$$
which, as above, implies 
$$\Ga(\beta)\,E_{\a,\beta}(-x)\; \le\; \frac{1}{1 +\frac{\Ga(\beta)}{\Ga(\beta+\a)}\, x}$$
for every $x\ge 0.$ 

The argument for the other inequality is analogous to that of Proposition \ref{KSO3}. By density, we only need to consider the case $\a = p/n$ and $\beta = (p+q)/n$ with $p<n$ and $q$ positive integers. By (\ref{MomXab}) and the Legendre-Gauss multiplication formula, we obtain 
$$\esp \lcr (\Y_{\a,1,l})^{ns}\rcr\; =\;\frac{n^{ns}}{p^{ps}}\,\times \,\esp\lcr \lpa\B_{\frac{1}{n},\frac{q}{np}}\rpa^s\rcr\,\times\,\frac{\prod_{i=2}^{n} (in^{-1})_s}{\prod_{j=1}^{p-1} (jp^{-1} + (p+q)(np)^{-1})_s}$$
for every $s > -n^{-1}.$ On the other hand, one has
$$\esp\lcr \lpa \frac{\Ga(\beta-\a)}{\Ga(\beta)}\;\L\rpa^{\! ns}\,\rcr \; = \; \frac{n^{ns}}{p^{ps}}\;\esp\lcr \lpa\kappa_{\a,\beta}\,\times \,\GG_{\frac{1}{n}}\,\times\, \prod_{j=1}^{p-1}\; \GG_{\frac{j}{p} +\frac{p+q}{np}}\rpa^{\!\! s}\,\rcr\,\times\,\frac{\prod_{i=2}^{n} (in^{-1})_s}{\prod_{j=1}^{p-1} (jp^{-1} + (p+q)(np)^{-1})_s}$$
with
$$\kappa_{\a,\beta}\; =\; p^p \lpa\frac{\Ga(qn^{-1})}{\Ga((p+q)n^{-1})} \rpa^n.$$
Comparing these two formulas, we are reduced to show 
$$(\B_{\frac{1}{n},\frac{q}{np}})^{\frac{1}{n}}\,\prec_{st}\, p^{\frac{p}{n}} \lpa\frac{\Ga(qn^{-1})}{\Ga((p+q)n^{-1})} \rpa\times\lpa\GG_{\frac{1}{n}}\,\times\, \prod_{j=1}^{p-1}\; \GG_{\frac{j}{p} +\frac{p+q}{np}}\rpa^{\frac{1}{n}}$$
for every $p< n$ and $q$ positive integers. This is obtained in the same way as above via the single intersection property. We leave the details to the reader.

\endproof
 
\section{Asymptotic behaviour of fractional extreme densities} 

\label{ABTD}

In this section which is a complement to \cite{BSV}, we study the behaviour of the density functions of the fractional Weibull and Fr\'echet distributions at both ends of their support. To this end, we also evaluate their Mellin transforms in terms of Barnes' double Gamma function. Along the way, we give the exact asymptotics of $x\mapsto E_{\a,m,l} (x)$ on the negative half-line, in the completely monotonic case $\a\in [0,1]$ and $ l \ge m -1/\a.$

\subsection{The fractional Weibull case} 

In \cite{BSV}, a fractional Weibull distribution function with parameters $\a\in [0,1]$ and $\lbd, \rho > 0$ is defined as the unique distribution function $F^\W_{\a,\lbd, \rho}$ on $(0,\infty)$ solving the fractional differential equation
$$\Dap F (x)\; =\; \lbd\, x^{\rho -\a} \bF (x)$$
where $\bF = 1-F$ denotes the associated survival function and $\Dap$ a progressive Liouville fractional derivative on $(0,\infty).$ The case $\a = 1$ corresponds to the standard Weibull distribution. In \cite{BSV}, it is shown that this distribution function exists and is given by
$$F^\W_{\a,\lbd, \rho} (x) \; =\; 1\; -\; E_{\alpha,\frac{\rho}{\a},\frac{\rho}{\a}-1}(-\lambda x^\rho)$$
for every $x \ge 0$ - see the formula following (3.1) in \cite{BSV}. In particular, the density $f^\W_{\a,\lbd, \rho}$ is real-analytic on $(0,\infty)$ and has the following asymptotic behaviour at zero: 
$$f^\W_{\a,\lbd, \rho}(x) \; \sim\; \lpa\frac{\lbd\,\Ga(\rho +1-\a)}{\Ga(\rho)}\rpa x^{\rho -1}\quad\mbox{as $x\to 0.$}$$
The behaviour of $f^\W_{\a,\lbd, \rho}$ at infinity is however less immediate, and to this aim we will need an exact expression for the Mellin transform of the random variable $\Warl$ with distribution function $F^\W_{\a,\lbd, \rho},$ which has an interest in itself.

\begin{prop}
\label{MASW}
The Mellin transform of $\Warl$ is 
$$\esp\lcr\Warl^s\rcr\; =\; \lpa\frac{\rho^\a}{\lbd}\rpa^{\frac{s}{\rho}}\Ga(1+ s\rho^{-1})\,\times\,\frac{[\rho +(1-\a);\rho]_{-s}}{[\rho;\rho]_{-s}}$$
for every $s\in (-\rho,\rho).$ As a consequence, one has
$$f^\W_{\a,\lbd, \rho}(x)\; \sim\; \lpa\frac{\rho }{\lbd \Ga(1-\a)}\rpa x^{-\rho -1} \quad \mbox{as $x\to\infty.$}$$
\end{prop}

\proof
We start with a more concise expression of (\ref{MellY}) for $l = m-1,$ which is a direct consequence of (\ref{Conca4}):
$$\esp[(\Y_{\a,\frac{\rho}{\a},\frac{\rho}{\a} -1})^s]\; =\; \rho^{-s} \times\frac{[1+(1-\a)\rho^{-1};\rho^{-1}]_s}{[1;\rho^{-1}]_s}\cdot$$
By Theorem 1.1 in \cite{BSV} and using the notations therein, we deduce
\begin{eqnarray*}
\esp\lcr\Warl^s\rcr & = &\esp\lcr\lpa \frac{\L}{\lbd \Y_{\a,\frac{\rho}{\a},\frac{\rho}{\a} -1}}\rpa^{\frac{s}{\rho}}\rcr\\
& = & \lpa\frac{\rho}{\lbd}\rpa^{\frac{s}{\rho}}\Ga(1+ s\rho^{-1})\,\times\, \frac{[1+(1-\a)\rho^{-1};\rho^{-1}]_{-s\rho^{-1}}}{[1;\rho^{-1}]_{-s\rho^{-1}}}\\
& = & \lpa\frac{\rho^\a}{\lbd}\rpa^{\frac{s}{\rho}}\Ga(1+ s\rho^{-1})\,\times\,\frac{[\rho +(1-\a);\rho]_{-s}}{[\rho;\rho]_{-s}}
\end{eqnarray*}
for every $s\in (-\rho,\rho)$ as required, where the third equality comes from (\ref{Conca3}). The asymptotic behaviour of the density at infinity is then a standard consequence of Mellin inversion. First, we observe from the above formula and (\ref{zeroes}) that the first positive pole of $s\mapsto \esp\lcr\Warl^s\rcr$ is simple and isolated in the complex plane at $s = \rho,$ with
\begin{eqnarray*}
\esp\lcr\W_{\a,\lbd,\rho}^s\rcr & \sim & \lpa\frac{\rho^\a}{\lbd}\rpa\,\times\,\frac{[\rho +(1-\a);\rho]_{-\rho}}{[\rho;\rho]_{-s}} \\
& \sim & \lpa\frac{\rho^{\rho +\a}}{\lbd}\rpa\,\times\,\frac{[\rho +(1-\a);\rho]_{-\rho}}{[2\rho;\rho]_{-\rho}}\,\times\,(\rho)_{-s}\; = \; \frac{\rho\, \Gamma (\rho-s)}{\lbd \Ga(1-\a)}\;\sim\; \frac{\rho}{\lbd \Ga(1-\a)\, (\rho -s)}
\end{eqnarray*}
as $s\uparrow \rho,$ where the second asymptotics comes from (\ref{Conca4}) and the equality from (\ref{Conca2}). Therefore, applying Theorem 4 (ii) in \cite{FGD1995} - beware the correction $(\log x)^k \to (\log x)^{k-1}$ to be made in the expansion of $f(x)$ therein, we obtain
$$f^\W_{\a,\lbd,\rho}(x)\; \sim\; \lpa\frac{\rho }{\lbd \Ga(1-\a)}\rpa x^{-\rho -1} \quad \mbox{as $x\to\infty$}$$
as required.
\endproof

\begin{rem}
\label{rW2}
{\em (a) Another proof of the asymptotic behaviour at infinity can be obtained from that of the so-called generalized stable densities. More precisely, using the identity in law on top of p.12 in \cite{BSV} and the notation therein, we see by multiplicative convolution, having set $f^\cG_{\a,\rho}$ for the density of the generalized stable random variable $\cG(\rho +1-\a,1-\a),$ that
\begin{eqnarray*}
f^\W_{\a,\lbd,\rho}(x) & = & \lbd \, x^{\rho-1} \int_0^\infty f^\cG_{\a,\rho}(y)\,y^{-\rho}\, e^{-\frac{\lbd}{\rho}(\frac{x}{y})^\rho} dy \\
& = & \lpa\frac{\lbd}{\rho}\rpa^{\frac{1}{\rho}}  \int_0^\infty f^\cG_{\a,\rho}\lpa x(\rho \lbd^{-1} t)^{-\frac{1}{\rho}}\rpa\,t^{-\frac{1}{\rho}}\,e^{-t}\, dt\\
& \sim & \lpa\frac{\rho }{\lbd \Ga(1-\a)}\int_0^\infty t \,e^{-t}\, dt \rpa x^{-\rho -1}\; =\; \lpa\frac{\rho }{\lbd \Ga(1-\a)}\rpa x^{-\rho -1}
\end{eqnarray*}
as $x\to\infty,$ where for the asymptotics we have used the Proposition in \cite{JSW2018} and a direct integration. This argument does not make use of Mellin inversion and is overall simpler than the above. However, it does not convey to the fractional Fr\'echet case. 

\medskip

(b) The Mellin transform simplifies for $\a = 0$ and $\a = 1\! :$ using (\ref{Conca1}) and (\ref{Rhorho}) we recover
$$\esp[\W_{0,\lbd,\rho}^s] \; = \; \lbd^{-\frac{s}{\rho}}\Ga(1+s\rho^{-1})\Ga(1-s\rho^{-1})\qquad\mbox{and}\qquad \esp[\W_{1,\lbd,\rho}^s] \; = \; \lpa\frac{\rho}{\lbd}\rpa^{\frac{s}{\rho}}\Ga(1+s\rho^{-1})$$
in accordance with the scaling property $\Warl \elaw \lbd^{-1/\rho}\W_{\a,1,\rho}$ and the identities given at the bottom of p.3 in \cite{BSV}. The Mellin transform takes a simpler form in two other situations. 

\begin{itemize}

\item For $\rho=\a,$ we obtain from (\ref{MellY}), (\ref{Conca1}) and (\ref{Conca2})
$$\esp[(\Y_{\a,1,0})^s]\; =\; \frac{\Ga(1+s)}{\Ga(1+\a s)}\; =\; \esp[\Z_\a^{-\a s}],$$
in accordance with Remark 3.1 (d) in \cite{BSV}. This yields
$$\W_{\a,\lbd,\a}\;\elaw\; \lpa \frac{\L}{\lbd \Y_{\a,1,0}}\rpa^{\frac{1}{\a}}\;\elaw\; \lbd^{-\frac{1}{\a}}\,\Z_\a\,\times\,\L^{\frac{1}{\a}},$$ 
an identity which was already discussed for $\lbd = 1$ in the introduction of \cite{BSV} as the solution to (1.3) therein. The Mellin transform reads
$$\esp[\W_{\a,\lbd,\a}^s] \; = \; \lbd^{-\frac{s}{\a}}\,\frac{\Ga(1+s\a^{-1})\Ga(1-s\a^{-1})}{\Ga(1-s)}\cdot$$

\item For $\rho=1-\a,$  where we obtain from (\ref{Conca2}) 
$$\esp[\W_{1-\rho,\lbd,\rho}^s] \; = \; \lpa\frac{\rho}{\lbd}\rpa^{\frac{s}{\rho}}\,\frac{\Ga(1+s\rho^{-1})\Ga(\rho -s)}{\Ga(\rho)}\qquad\mbox{and}\qquad \W_{1-\rho,\lbd,\rho}\; \elaw\; \lpa\frac{\rho}{\lbd}\rpa^{\frac{1}{\rho}}\,\L^{\frac{1}{\rho}}\,\times\,\GG^{-1}_{\rho}.$$

\end{itemize}

\medskip

(c) The two cases $\rho = \a$ and $\rho = 1-\a$ have a Mellin transform expressed as the quotient of a finite number of Gamma functions. This makes it possible to use a Mellin-Barnes representation of the density in order to get its full asymptotic expansion at infinity. Using the standard notation of Definition C.1.1 in \cite{AAR1999}, one obtains 
$$f^\W_{\a,\lbd,\a}(x)\; \sim\; \sum_{n\ge 1}\, \frac{n\a\, x^{-1-n\a}}{\lbd^n\Ga(1-n\a)}\qquad\mbox{and}\qquad f^\W_{1-\a,\lbd,\a}(x)\; \sim\; \frac{\a x^{-\a -1}}{\lbd \Ga(\a)} \sum_{n\ge 0}\, (-1)^n\, \frac{\Ga\lpa\frac{n}{\rho} +2\rpa}{n!}\!\lpa\frac{\lbd}{\rho}\rpa^{-\frac{n}{\rho}}\!\!\! x^{-n}$$
which are everywhere divergent. The first expansion can also be obtained from (1.8.28) in \cite{KST} using
$$f^\W_{\a,\lbd,\a}(x)\; =\; \lbd\, x^{\a-1} \,E_{\a,\a}(-\lbd x^\a).$$
Unfortunately, the Mellin transform of $\Warl$ might have poles of variable order and it seems difficult to obtain a general formula for the full asymptotic expansion at infinity of $f^\W_{\a,\lbd,\rho}(x)$.}
\end{rem}

Writing
$$E_{\a,\frac{\rho}{\a}, \frac{\rho}{\a}-1}(-\lbd x^\rho) \; =\; \pb[\Warl > x] \; =\; \int_x^\infty f^\W_{\a,\lbd,\rho}(y)\, dy,$$ 
we obtain by integration the following asymptotic behaviour at infinity, which is valid for any $\a\in (0,1]$ and $m > 0:$ 
$$E_{\a,m, m-1}(-x)\;\sim\; \frac{1}{\Ga(1-\a)\, x}\qquad\mbox{as $x\to\infty.$}$$ 
This behaviour, which turns out to be the same as that of the classical Mittag-Leffler function $E_\a(-x)$ - see e.g. (3.4.15) in \cite{GKMR}, gives the reason why the constant in the lower bound of Theorem \ref{KSO1} is optimal - see the above Remark \ref{rKS1} (a). It is actually possible to get the exact behaviour of $E_{\a,m, l}(-x)$ at infinity for any $\a\in (0,1], m > 0$ and $l> m - 1/\a.$ We include this result here since it seems unnoticed in the literature on Kilbas-Saigo functions.

\begin{prop}
\label{ASKS}
For any $\a\in [0,1], m > 0$ and $l> m - 1/\a,$ one has
$$E_{\a,m, l}(-x)\;\sim\; \frac{\Ga(1+\a(l+1-m))}{\Ga(1+\a(l-m))\, x}\qquad\mbox{as $x\to\infty.$}$$ 
\end{prop}

\proof 

The case $\a = 0$ is obvious since $E_{0,m,l} (x) = 1/(1-x).$ For $\a \in (0,1],$ setting $\delta = 1/\a m,$ recall from (\ref{Marr}) that for every $s \in (0,1)$ one has 
\begin{eqnarray*}
\int_0^\infty E_{\a, m, l}(-x)\, x^{s-1}\, dx & = & \Ga(s)\Ga(1-s)\, \times\,\frac{[(\a l +1)\delta;\delta]_{-s}}{[1/m +(\a l+1)\delta;\delta]_{-s}} \\
& \sim & \frac{[(\a l +1)\delta;\delta]_{-1}}{[1/m +(\a l+1)\delta;\delta]_{-1}\, (1-s)}\; = \; \frac{\Ga(1+\a(l+1-m))}{\Ga(1+\a(l-m))\, (1-s)}
\end{eqnarray*}
as $s\uparrow 1,$ where in the equality we have used the concatenation formula (\ref{Conca1}). The asymptotic behaviour follows then by Mellin inversion as in the proof of Proposition \ref{MASW}.

\endproof

\begin{rem}
\label{LimitF}
{\em In the boundary case $l = m-1/\a,$ the behaviour of $E_{\a, m, m-1/\a}(-x)$ at infinity, which has different speed and a more complicated constant, will be obtained with the help of the fractional Fr\'echet distribution - see Remark \ref{rF2} (c) below.} 
\end{rem}
\endproof

We end this paragraph with the following conjecture which is natural in view of Proposition \ref{ASKS}. We know by Theorem \ref{KSO1} resp. Proposition \ref{KSO4} that this conjecture is true for the cases $l = m -1$ and $m = 1.$

\begin{conj}
\label{OrdKS}
For every $\a\in(0,1], m > 0, l > m-1/\a$ and $x\ge 0,$ one has 
$$\frac{1}{1 +\frac{\Ga(1+\a(l-m))}{\Ga(1+\a(l-m+1))}\, x}\;\le\; E_{\a,m,l}(-x)\; \le\; \frac{1}{1 +\frac{\Ga(1+\a l)}{\Ga(1+\a(1+l))}\, x}\cdot$$
\end{conj}

\medskip

\subsection{The Fr\'echet case} In \cite{BSV}, a fractional Fr\'echet distribution function with parameters $\a\in [0,1]$ and $\lbd, \rho > 0$ is defined as the unique distribution function $F^\F_{\a,\lbd, \rho}$ on $(0,\infty)$ solving the fractional differential equation
$$\Dm \bF(x)\; =\; \lbd\, x^{-\rho -\a} F(x)$$
where $\Dm$ denotes a regressive Liouville fractional derivative on $(0,\infty).$ The case $\a = 1$ corresponds to the standard Fr\'echet distribution. In \cite{BSV}, it is shown that this distribution function exists and is given by
$$F^\F_{\a,\lbd, \rho} (x) \; =\; E_{\alpha,\frac{\rho}{\a},\frac{\rho-1}{\a}}(-\lambda x^{-\rho})$$
for every $x \ge 0$ - see the formula following (3.4) in \cite{BSV}. In particular, the density $f^\F_{\a,\lbd, \rho}$ is real-analytic on $(0,\infty)$ and has the following asymptotic behaviour at infinity: 
$$f^\F_{\a,\lbd, \rho}(x) \; \sim\; \lpa\frac{\lbd\,\Ga(\rho +1)}{\Ga(\rho+\a)}\rpa x^{-\rho -1}\quad\mbox{as $x\to \infty.$}$$
The behaviour of the density at zero is less immediate and we will need, as in the above paragraph, the exact expression of the Mellin transform of the random variable $\Farl$ with distribution function $F^\F_{\a,\lbd, \rho},$ whose strip of analyticity is larger than that of $\Warl.$ 

\begin{prop}
\label{MASF}
The Mellin transform of $\Farl$ is 
$$\esp\lcr\Farl^s\rcr\; =\; \lpa\frac{\rho^\a}{\lbd}\rpa^{-\frac{s}{\rho}}\Ga(1- s\rho^{-1})\,\times\,\frac{[\rho +1;\rho]_{s}}{[\rho +\a;\rho]_{s}}$$
for every $s\in (-\rho-\a,\rho).$ As a consequence, one has
$$f^\F_{\a,\lbd, \rho}(x)\; \sim\; \lpa\frac{\rho^{\frac{\a^2}{\rho}}(\rho +\a) }{\lbd^{1 +\frac{\a}{\rho}}}\,\Ga(1+\a)\, G(1-\a; \rho)\, G(1+\a; \rho)\rpa x^{\rho +\a -1} \quad \mbox{as $x\to 0.$}$$
\end{prop}

\proof The evaluation of the Mellin transform is done as for the fractional Weibull distribution, starting from the expression
$$\esp[(\Y_{\a,\frac{\rho}{\a},\frac{\rho -1}{\a}})^s]\; =\; \rho^{-s} \times\frac{[1+\rho^{-1};\rho^{-1}]_s}{[1 +\a\rho^{-1};\rho^{-1}]_s}$$
which is a consequence of (\ref{MellY}) and (\ref{Conca2}). By Theorem 1.2 in \cite{BSV} and (\ref{Conca3}), we obtain the required formula
$$\esp\lcr\Farl^s\rcr \; = \; \esp\lcr\lpa \frac{\L}{\lbd \Y_{\a,\frac{\rho}{\a},\frac{\rho-1}{\a}}}\rpa^{-\frac{s}{\rho}}\rcr\; =\;\lpa\frac{\rho^\a}{\lbd}\rpa^{-\frac{s}{\rho}}\Ga(1- s\rho^{-1})\,\times\,\frac{[\rho +1;\rho]_{s}}{[\rho +\a;\rho]_{s}}\cdot$$
The asymptotic behaviour of $f^\F_{\a,\lbd, \rho}(x)$ at zero follows then as for that of $f^\W_{\a,\lbd, \rho}(x)$ at infinity, in considering the residue at the first negative pole $s = -(\rho+\a)$ which is simple and isolated in the complex plane, applying Theorem 4 (i) in \cite{FGD1995} - with the same correction as above, and making various simplifications. We omit details. 

\endproof

\begin{rem}
\label{rF2}{\em (a) Comparing the Mellin transforms, Propositions \ref{MASW} and \ref{MASF} imply the factorization
\begin{equation}
\label{FacWF}
\Warl^{-1}\;\elaw\; \Farl\;\times\;\Z(\rho +1-\a,\rho+\a;\rho,\rho+1;\rho).
\end{equation}
In general, it follows from Theorem \ref{KSCM} that for every $\a\in (0,1], m, \lbd > 0$ and $l > m- 1/\a,$ there exists a positive random variable having distribution function $E_{\a,m, l}(-\lbd x^{-\a m}),$ and which is given by (\ref{MellY}), (\ref{BKS}) and Theorem 1.2 in \cite{BSV} as the independent product
$$\F_{\a,\lbd, \a m}\;\times\;(\X_{\a,m,l})^{\frac{1}{\a m}} \;\elaw\; \F_{\a,\lbd, \a m}\;\times\;\Z(\a l +1,\a (m+1);\a m,\a(l+1) +1;\a m),$$
where the identity in law follows from (\ref{Conca3}). In this respect, the fractional Fr\'echet distributions can be viewed as the ``ground state'' distributions associated to the Kilbas-Saigo functions $E_{\a,m, l},$ in the boundary case $l = m-1/\a.$ 

\medskip

(b) As above, the Mellin transform simplifies for $\a = 0,1\! :$ we get
$$\esp[\F_{0,\lbd,\rho}^s] \; = \; \lbd^{\frac{s}{\rho}}\,\Ga(1+s\rho^{-1})\Ga(1-s\rho^{-1})\qquad\mbox{and}\qquad \esp[\F_{1,\lbd,\rho}^s] \; = \; \lpa\frac{\lbd}{\rho}\rpa^{\frac{s}{\rho}}\Ga(1-s\rho^{-1}),$$
in accordance with the scaling property $\Farl \elaw \lbd^{1/\rho}\F_{\a,1,\rho}$ and the identities given after the statement of Theorem 1.2 in \cite{BSV}.
The Mellin transform also takes a simpler form in the same other situations as above.

\begin{itemize}

\item For $\rho=\a,$ with
$$\esp[\F_{\a,\lbd,\a}^s] \; = \; \lbd^{\frac{s}{\a}}\,\frac{\Ga(\a)\Ga(1+s\a^{-1})\Ga(1-s\a^{-1})}{\Ga(\a+s)}\cdot$$
This yields the identity $\F_{\a,\lbd,\a}\,\elaw\, \lbd^{\frac{1}{\a}}\,(\Z^{-1}_\a)^{(\a)}\,\times\,\L^{-\frac{1}{\a}},$ which was discussed for $\lbd =1$ in the introduction of \cite{BSV} as the solution to (1.4) therein. This is also in accordance with Remark 3.3 (c) in \cite{BSV}, since 
$$(T_\a^{(1)})^{\frac{1}{\a}}\;\elaw\; ((\Z^{-\a}_\a)^{(1)})^{\frac{1}{\a}}\;\elaw\;(\Z^{-1}_\a)^{(\a)}.$$ 
Notice that the constant appearing in the asymptotic behaviour of the density at zero is also simpler: one finds
\begin{equation}
\label{Fala}
f^\F_{\a,\lbd, \a}(x)\; \sim\; \lpa\frac{2\a\,\Ga(1+\a)}{\lbd^2\,\Ga(1-\a)}\rpa x^{2\a -1} \quad \mbox{as $x\to 0.$}
\end{equation}

\item For $\rho=1-\a,$ with
$$\esp[\F_{1-\rho,\lbd,\rho}^s] \; = \; \lpa\frac{\lbd}{\rho}\rpa^{\frac{s}{\rho}}\,\Ga(1-s \rho^{-1})\Ga(1+s)\quad\mbox{and}\quad\F_{1-\rho,\lbd,\rho}\; \elaw\; \lpa\frac{\lbd}{\rho}\rpa^{\frac{1}{\rho}}\,\L^{-\frac{1}{\rho}}\,\times\,\L.$$
Here, the density converges at zero to a simple constant: one finds
$$f^\F_{1-\rho,\lbd, \rho}(x)\; \to\; \lpa\frac{\rho}{\lbd}\rpa^{\frac{1}{\rho}} \Ga(1+\rho^{-1})\qquad \mbox{as $x\to 0.$}$$
\end{itemize}

\medskip

(c) Integrating the density and using $\pb[\Farl \le x] = E_{\a,\frac{\rho}{\a}, \frac{\rho -1}{\a}}(-\lbd x^{-\rho}),$ we obtain the following asymptotic behaviour at infinity for any $\a\in (0,1]$ and $m > 0,$ which is more involved than that of Proposition \ref{ASKS}: 
$$E_{\a,m, m-\frac{1}{\a}}(-x)\;\sim\; (\a m)^{\frac{\a}{m}}\Ga(1+\a)\, G(1-\a; \a m)\, G(1+\a; \a m) \, x^{-1-\frac{1}{m}}\qquad\mbox{as $x\to\infty.$}$$ 
For $m=1,$ this behaviour matches the first term in the full asymptotic expansion
$$E_{\a,1, 1-\frac{1}{\a}}(-x)\; =\; \Ga(\a) E_{\a,\a}(-x)\;\sim\; \Ga(\a)\sum_{n\ge 1} \frac{(-1)^n}{\Ga(-\a n)\, x^{n+1}}\cdot$$
As for $E_{\a,m, m-1}(-x),$ a full asymptotic expansion of $E_{\a,m, m-\frac{1}{\a}}(-x)$ at infinity seems difficult to obtain for all values of $m.$}
\end{rem}

\section{Some complements on the Le Roy function} 

\label{SCLR}

In this section we show some miscellaneous results on the Le Roy function
$$\cL_\a(x) \; =\; \sum_{n\ge 0} \frac{x^n}{(n!)^\a}, \qquad \a >0, \; x\in\rl.$$
In \cite{BSV}, this function played a role in the construction of a fractional Gumbel distribution - see Theorem 1.3 therein. The Le Roy function, which has been much less studied than the classical Mittag-Leffler function, can be viewed as an alternative generalization of the exponential function. Throughout, we giscard the explicit case $\cL_1(x) = E_1(x) = e^x.$\\

We begin with the asymptotic behaviour at infinity. Le Roy's original result - see \cite{L1899} p. 263 - reads 
$$\cL_\a (x)\; \sim\; \frac{(2\pi)^{\frac{1-\a}{2}}}{\sqrt{\a}}\, x^{\frac{1-\a}{2}}\, e^{\a x^{\frac{1}{\a}}}\qquad \mbox{as $x\to\infty,$}$$
and is obtained by a variation on Laplace's method. An extension of this asymptotic behaviour was recently given in \cite{G} for the so-called Mittag-Leffler functions of Le Roy type. Laplace's method can also be used to solve Exercise 8.8.4 in \cite{Olver}, which states
\begin{equation}
\label{Olv1}
\cL_\a (-x)\; = \; \frac{2 (2\pi)^{\frac{1-\a}{2}}}{\sqrt{\a}}\,  x^{\frac{1-\a}{2\a}}\, e^{\a \cos(\pi/\a) x^{\frac{1}{\a}}}\, \lpa \sin\lpa \pi/2\a + \a \sin\lpa\pi/\a\rpa x^{\frac{1}{\a}}\rpa\, +\, O(x^{-\frac{1}{\a}})\rpa
\end{equation}
for $\a\ge 2$ and 
\begin{equation}
\label{Olv2}
\cL_\a (-x)\; \sim \; \frac{1}{\a^\a\, \Ga(1-\a) \, x\, (\log x)^\a}
\end{equation}
for $\a \in (1,2),$ as $x\to \infty.$ The following estimate, which seems to have passed unnoticed in the literature, completes the picture.

\begin{prop}
\label{Olv3}
For every $\a\in (0,1),$ one has
$$\cL_\a (-x)\; \sim \; \frac{1}{\Ga(1-\a) \, x\, (\log x)^\a}\qquad \mbox{as $x\to\infty.$}$$
\end{prop}

\proof

In the proof of Theorem 1.3 in \cite{BSV} it is shown that
$$\cL_\a (-x)\; =\; \pb\lcr \L > x \L_\a\rcr\; =\; \int_0^\infty e^{-xt} \, f_\a(t)\, dt$$
where 
$$\L_\a\; \elaw\; \int_0^\infty e^{-\sga_t}\, dt$$ 
has density $f_\a$ on $(0,\infty)$ and Mellin transform 
$$\esp \lcr \L_\a^s\rcr\; =\; \Ga(1+s)^{1-\a}, \qquad s > -1.$$
In particular, we have $f_\a = e_{1-\a}$ with the notation of \cite{BL2015} and Theorem 2.4 therein implies
\begin{equation}
\label{Asfa}
f_\a(x)\;\sim\; \frac{1}{\Ga(1-\a)\, (-\log x)^\a}\qquad \mbox{as $x\to 0.$}
\end{equation}
Plugging this estimate into the above expression for $\cL_\a (-x),$ we conclude the proof by a direct integration.

\endproof

\begin{rem}
\label{Rolv}
{\em (a) The estimate (\ref{Asfa}) also gives the asymptotic behaviour, at the right end of the support, of the density of the fractional Gumbel random variable $\Gal$ which is defined in Theorem 1.3 of \cite{BSV}. Indeed, by the definition and multiplicative convolution the density of $e^{\lbd \Gal}$ on $(0,\infty)$ writes
$$\int_0^\infty e^{-xy} \; y f_\a(y)\, dy \; \sim\; \frac{1}{\Ga(1-\a)\, x^2\, (\log x)^\a}\qquad\mbox{as $x\to\infty,$}$$
where the estimate follows from (\ref{Asfa}) as in the proof of Proposition \ref{Olv3}. A change of variable implies then
$$f^\G_{\a,\lbd}(x)\;\sim\;\lpa\frac{\lbd^{1-\a}}{\Ga(1-\a)}\rpa\, x^{-\a}\, e^{-\lbd x}\qquad\mbox{as $x\to\infty.$}$$
Notice that at the left end of the support, there is a convergent series representation which is given by Corollary 3.6 in \cite{BSV}. 

\medskip
 
(b) In the case $\a =2,$ one has $\cL_2(x) = I_0(2\sqrt{x})$ and $\cL_2(-x) = J_0(2\sqrt{x})$ for all $x\ge 0,$ where $I_0$ and $J_0$ are the classical, modified or not, Bessel functions with index 0. In particular, a full asymptotic expansion for $\cL_2$ at both ends of the support is available, to be deduced e.g. from (4.8.5) and (4.12.7) in \cite{AAR1999}. These expansions also exist when $\a$ is an integer since $\cL_\a$ is then a generalized Wright function - see Chapter F.2.3 in \cite{GKMR} and the original articles by Wright quoted therein. The case when $\a$ is not an integer does not seem to have been investigated, and might be technical in the absence of a true Mellin-Barnes representation.}
\end{rem}

Our next result characterizes the connection between the entire function $\cL_\a(z)$ and random variables. Recall that a function $f : \CC\to\CC$ which is holomorphic in a neighbourhood $\Omega$ of the origin is a moment generating function (MGF) if there exists a real random variable $X$ such that 
$$f(z) \; =\; \esp\lcr e^{z X}\rcr,\qquad z\in\Omega.$$
In particular, it is clear that $\cL_0$ is the MGF of the exponential law $\L$ and $\cL_1$ that of the constant variable $\Un.$ The following provides a characterization.

\begin{prop}
\label{LRCM}
The function $\cL_\a(z)$ is the {\em MGF} of a real random variable if and only if $\a\le 1.$ In this case, one has
$$\cL_\a(z)\; =\; \esp\lcr e^{z \L_\a}\rcr, \qquad z\in\CC.$$
\end{prop}

\proof

The if part is a direct consequence of the proof of Proposition \ref{Olv3}. On the other hand, the estimates (\ref{Olv1}) and (\ref{Olv2}) show that $\cL_\a(z)$ takes negative values on $\rl^-,$  so that it cannot be the moment generating function of a real random variable, when $\a  > 1.$ This completes the proof.

\endproof

Observe that since $\L_\a$ is non-negative, the above result also shows $\cL_\a(-x)$ is CM on $(0,\infty)$ if and only if $\a\le 1,$ echoing Pollard's aforementioned classical result for the Mittag-Leffler $E_\a(-x).$ One can ask whether there are further complete monotonicity properties for $\cL_\a,$ as in \cite{TS2015} for $E_\a.$ Our last result for the Le Roy function is a monotonicity property which is akin to Proposition \ref{Mono}.

\begin{prop}
\label{LRMono}
The mapping $\a\mapsto\cL_\a (x)$ is non-increasing on $[0,1]$ for every $x\in\rl.$ 
\end{prop}
 
\proof The fact that $\a\mapsto\cL_\a (x)$ decreases on $\rl^+$ is obvious for $x\ge 0,$ by the definition of $\cL_\a.$ To show the property on $[0,1]$ for $x <0,$ we will use a convex ordering argument. More precisely, the Malmsten formula (\ref{Malm1}) and the L\'evy-Khintchine formula show that for every $t \in [0,1],$ the random variable $\G_{1-t} = \log \L_{1-t}$ is the marginal at time $t$ of a real L\'evy process, since $\esp[e^{\ii z \G_{1-t}}] =  \Ga(1+ \ii z)^t = e^{t \psi(z)}$ for every $z\in\rl,$ with
$$\psi(z)\; =\; -\gamma\ii z\, +\, \int^0_{-\infty} (e^{\ii zx} - 1 - \ii zx)\, \frac{dx}{\vert x\vert (e^{\vert x\vert} -1)}\cdot$$  
This is actually well-known - see Example E in \cite{CPY}. By independence and stationarity of the increments of a L\'evy process, we deduce that there exists a multiplicative martingale $\{M_t, t\in [0,1]\}$ such that $M_t \elaw \L_{1-t}$ for every $t \in [0,1]$ and Jensen's inequality implies
$$\L_\beta\;\prec_{cx}\; \L_\a$$
for every $0\le\a\le\beta\le 1.$ Applying the definition of convex ordering to the function $\varphi(x) = e^x,$ we get $\cL_\beta(x)\le\cL_\a(x)$ for every $x<0$ and $0\le\a\le\beta\le 1,$ as required. 

\endproof

\begin{rem}
\label{Picq}

{\em (a) In the terminology of \cite{HPRY}, the family $\{\L_{1-\a}, \; \a\in[0,1]\}$ is a peacock, whose associated multiplicative martingale is completely explicit. We refer to \cite{HPRY} for numerous examples of explicit peacocks related to exponential functionals of L\'evy processes. Observe from Lemma \ref{Tcx} that the family $\{\T(a,b,t), \; t >0\}$ is also a peacock.

\medskip

(b) Letting $\a\to 0$ and $\a\to 1$ in Proposition \ref{LRMono} leads to the bounds
$$e^{x}\;\le\;\cL_\beta(x)\;\le\;\cL_\a(x)\;\le\; \frac{1}{(1-x)_+}$$
for every $x \in \rl$ and  $0 < \a < \beta < 1.$ The hyperbolic upper bound is optimal as in Propositions \ref{KSO1} and \ref{KSO2}, because $\cL_\a(x)- 1\sim x$ as $x\to 0.$ The exponential lower bound is thinner than the order given in Proposition \ref{Olv3}. On the other hand, it does not seem that stochastic ordering arguments can help for a uniform estimate involving a logarithmic term.}
\end{rem}

It is natural to ask if the statement of Proposition \ref{LRMono} is also true for the classical Mittag-Leffler function, and this problem seems still open. 

\begin{conj}
\label{MLMono}
The mapping $\a\mapsto\Ea (x)$ is non-increasing on $[0,1]$ for every $x\in\rl.$ 
\end{conj}
Numerical simulations suggest a positive answer. For $x <0,$ see the GIF animation given on the top of the english Wikipedia page of the Mittag-Leffler function \cite{Wiki}. It is clear by the definition that $\a\mapsto\Ea(x)$ is non-increasing for every $x \ge 0$ on $[\a_0, \infty),$ where $1 +\a_0 = 1.46163...$ is the location of the minimum of the Gamma function on $(0,\infty).$ A direct consequence of Theorem B in \cite{TS2014} is also that
$$\a\;\mapsto \; \Ea(\Ga(1+\a) x)$$
is non-increasing on $[1/2,1]$ for every $x\in\rl.$ The constant $\Ga(1+\a)$ appears above because of the convex ordering argument used in \cite{TS2014}. It seems that other kinds of arguments are necessary to study the monotonicity of $\a\mapsto \Ea(x)$ on $[0,1].$ \\

We would like to finish this paper with the following related monotonicity result, which relies on a stochastic ordering argument, for the generalized Mittag-Leffler function. 
  
\begin{prop}
\label{MonoEab}
For every $\a\in [0,1]$ and $x\in\rl,$ the mapping 
$$\beta\; \mapsto\;\Ga(\beta) E_{\a,\beta}(x)$$ 
is non-increasing on $(\a,\infty)$ if $x > 0$ and non-decreasing on $(\a,\infty)$ if $x < 0.$ 
\end{prop}

\proof

By Remark 3.3 (c) in \cite{BSV}, we have the probabilistic representation
$$\Ga(\beta) E_{\a,\beta}(x)\; =\; \esp \lcr e^{x \,\B_{\a,\beta -\a}^\a \times\, T_\a^{(1)}}\rcr$$
for every $\a\in [0,1], \beta > \a$ and $x\in\rl.$ Reasoning as in Proposition \ref{KSO3}, we see by factorization that it suffices to show that
$$\beta\;\mapsto\; \B_{\a,\beta -\a}^\a$$
is non-increasing on $(\a,\infty)$ for the usual stochastic order. On the other hand, the density function of the random variable $\B_{\a,\beta -\a}^\a$ is
$$\frac{\Ga(\beta)}{\Ga(\a + 1)\Ga(\beta -\a)}\, \lpa 1-x^{\frac{1}{\a}}\rpa ^{\beta - \a -1}
$$
on $[0,1)$ and its value at zero is by the log-convexity of the Gamma function an increasing function of $\beta.$ Moreover, the density functions of $\B_{\a,\beta -\a}^\a$ and $\B_{\a,\beta' -\a}^\a$ cross only once for $\beta\neq\beta',$ at
$$1\; -\; \lpa\frac{\Ga(\beta)\Ga(\beta'-\a)}{\Ga(\beta')\Ga(\beta -\a)} \rpa^{\frac{1}{\beta'-\beta}}\; \in \; (0,1).$$
The single intersection property finishes then the argument, as for Proposition \ref{KSO3}.

\endproof

\appendix
\section*{Appendix}

\renewcommand{\theequation}{A\thesection.\arabic{equation}}

\setcounter{equation}{0}

In this Appendix we recall some properties of Barnes' double Gamma function $G(z;\delta),$ which are used throughout the paper. For every $\delta > 0,$ this function is defined as the unique solution to the functional equation
\begin{equation}
\label{Conca1}
G(z+1;\delta)\; =\; \Gamma(z\delta^{-1}) G(z;\delta)
\end{equation}
with normalization $G(1;\delta) =1.$ The function is holomorphic on $\CC$ and admits the following Malmsten type representation
\begin{equation}
\label{Malm}
G(z;\delta)\; =\; \exp\int_0^\infty \lpa \frac{1-e^{-zx}}{(1-e^{-x})(1-e^{-\delta x})} \, -\, \frac{ze^{-\delta x}}{1-e^{-\delta x}}\, +\, (z-1)(\frac{z}{2\delta} -1)e^{-\delta x}\, - \, 1\rpa\frac{dx}{x}
\end{equation}
which is valid for $\Re(z) >0$ - see (5.1) in \cite{BK1997}. Putting (\ref{Conca1}) and (\ref{Malm}) together and making some simplifications, we recover the standard Malmsten formula for the Gamma function
\begin{equation}
\label{Malm1}
\Ga(1+z)\; =\; \exp\lacc -\gamma z\, +\, \int^0_{-\infty} (e^{zx} - 1 - zx)\, \frac{dx}{\vert x\vert (e^{\vert x\vert} -1)}\racc
\end{equation}
for every $z >-1,$ where $\gamma$ is Euler's constant. The following Stirling type asymptotic behaviour 
\begin{equation}
\label{Stirling}
\log G(z;\delta)\; -\; \frac{1}{2 \delta} \lpa z^2 \log z\, -\, (\frac{3}{2} + \log \delta)\, z^2 \, -\,(1+\delta)\, z\log z\rpa\; -\; A \, z\; -\; B \log z\; \to\; C
\end{equation}
is valid for $\vert z\vert\to\infty$ with $\vert \arg (z) \vert < \pi,$ for some real constants $A, B$ and $C$ which are given in (4.5) of \cite{BK1997}. There is a second concatenation formula
\begin{equation}
\label{Conca2}
G(z+\delta;\delta)\; =\; (2\pi)^{(\delta-1)/2} \delta^{1/2 -z}\Gamma(z) G(z;\delta)
\end{equation}
which is valid for all $z\in\CC,$ the right-hand side being understood as an analytic extension when $z$ is a non-positive integer - see (4.6) in \cite{K2011} and the references therein. Observe that (\ref{Conca1}) and (\ref{Conca2}) lead readily to the closed formula
\begin{equation}
\label{Rhorho}
G(\delta;\delta)\; =\; G(1+\delta;\delta)\; =\; (2\pi)^{(\delta-1)/2} \delta^{-1/2}.
\end{equation}
In this paper we make an extensive use of the following Pochhammer type symbol
\begin{equation}
\label{Poch}
[a;\delta]_s\; =\; \frac{G(a +s;\delta)}{G(a;\delta)}
\end{equation}
which is well-defined for every $a,\delta > 0$ and $s > - a.$ The following formula
\begin{equation}
\label{Conca3}
[a\delta^{-1};\delta^{-1}]_{s\delta^{-1}} \; =\; (2\pi)^{s(1/\delta-1)/2}\, \delta^{s^2/2\delta - s(1+(1-2a)/\delta)/2}\,[a;\delta]_s
\end{equation}
can be deduced from (4.10) in \cite{K2011} - beware the different normalization for $G(1;\delta)$ therein which becomes irrelevant when considering the Pochhammer type symbol. Notice also that (\ref{Conca2}) yields
\begin{equation}
\label{Conca4}
\delta^s\,[a+\delta;\delta]_s\; =\; (a)_s \,
[a;\delta]_s
\end{equation}
with the standard notation
$$(a)_s \; =\; \frac{\Ga(a+s)}{\Ga(a)}$$
for the usual Pochhammer symbol. Finally, we observe from the double product representation of $G(z,\delta)$ - see e.g. (4.4) in \cite{K2011}, that for every $a, \delta > 0$ one has
\begin{equation}
\label{zeroes}
\inf\{s > 0, \;  [a;\delta]_{-s} = 0\}\; =\; a
\end{equation} 
and that this zero is simple and isolated on the complex plane.


\begin{thebibliography}{10}

\bibitem{AAR1999}
G.~E.~Andrews, R.~Askey and R.~Roy. {\em Special functions.} Cambridge University Press, Cambridge, 1999.

\bibitem{BL2015}
C.~Berg and J.~L.~L\'opez. Asymptotic behaviour of the Urbanik semigroup. {\em J. Approx. Theory} {\bf 195}, 109-121, 2015.

\bibitem{BK1997}
J.~Billingham and A.~C.~King. Uniform asymptotic expansions for the Barnes double gamma function. {\em Proc. Roy. Soc. London Ser. A} {\bf 453}, 1817-1829, 1997.

\bibitem{BSV}
L.~Boudabsa, T.~Simon and P.~Vallois. Fractional extreme distributions. {\em Elec. J. Probab.} {\bf 25}, No. 115, 1-20, 2020.

\bibitem{CPY}
P.~Carmona, F.~Petit and M.~Yor. On the distribution and asymptotic results for exponential functionals of L\'evy processes. In: {\em Exponential functionals and principal values related to Brownian motion}. Bibliot\'eca de la Revista Matem\'atica Iberoamericana, 73-121, 1997.

\bibitem{DMV}
E.~C.~de Oliveira, F.~Mainardi and J.~Vaz. Fractional models of anomalous relaxation based on the Kilbas and Saigo function. {\em Meccanica} {\bf 49} (9), 2049-2060, 2014.
 
\bibitem{D2010}
D.~Dufresne. $G$ distributions and the beta-gamma algebra. {\em Elec. J. Probab.} {\bf 15}, No. 71, 2163-2199, 2010.

\bibitem{FGD1995}
P.~Flajolet, X.~Gourdon and P.~Dumas. Mellin transforms and asymptotics: Harmonic sums. {\em Theoret. Comput. Sci.} {\bf 144}, 3-58, 1995.

\bibitem{G}
S.~Gerhold. Asymptotics for a variant of the Mittag-Leffler function. {\em Int. Transf. Spec. Funct.} {\bf 23} (6), 397-403, 2012.

\bibitem{GKMR}
R.~Gorenflo, A.~A.~Kilbas, F.~Mainardi and S.~V. Rogosin. {\em Mittag-Leffler functions, related topics and applications.} Springer Verlag, Heidelberg, 2014.

\bibitem{HPRY}
F.~Hirsch, C.~Profeta, B.~Roynette and M. Yor. {\em Peacocks and associated martingales, with explicit constructions.} Springer-Verlag, Mailand, 2011.

\bibitem{JSW2018}
W.~Jedidi, T.~Simon and M.~Wang. Density solutions to a class of integro-differential equations. {\em J. Math. Anal. Appl.} {\bf 458} (1), 134-152, 2018.

\bibitem{KS95}
A.~A.~Kilbas and M.~Saigo. On solution of integral equation of Abel-Volterra type. {\em Differ. Integral Equ.} {\bf 8} (5), 993-1011, 1995.

\bibitem{KST}
A.~A.~Kilbas, H.~M.~Srivastava and J.~J.~Trujillo. {\em Theory and applications of fractional differential equations.} North-Holland, Amsterdam, 2006.

\bibitem{K2011}
A.~Kuznetsov. On extrema of stable processes. {\em Ann. Probab.} {\bf 39} (3), 1027-1060, 2011.

\bibitem{L1899}
E.~Le Roy. Valeurs asymptotiques de certaines s\'eries proc\'edant suivant les puissances enti\`eres et positives d'une variable r\'eelle. {\em Darboux Bull.} {\bf 24} (2), 245-268, 1899.

\bibitem{LS}
J.~Letemplier and T.~Simon. On the law of homogeneous stable functionals. {\em ESAIM P \& S} {\bf 23}, 82-111, 2019.

\bibitem{M}
F.~Mainardi. On some properties of the Mittag-Leffler function $\Ea(-t^\a),$ completely monotone for $t >0$ with $0<\a <1.$ {\em Discrete Cont. Dyn. Syst. Ser. B} {\bf 19} (7), 2267-2278, 2014.

\bibitem{Olver}
F.~W.~J. Olver. {\em Asymptotics and special functions.} Academic Press, New York, 1974.
 
\bibitem{S1999}
K.~Sato. {\em L\'evy processes and infinitely divisible distributions.} Cambridge University Press, Cambridge, 1999. 

\bibitem{SSh}
M.~Shaked and J.~G.~Shanthikumar. {\em Stochastic orders and their applications.} Springer Verlag, New York, 2007.

\bibitem{TS2014}
T.~Simon. Comparing Fr\'echet and positive stable laws. {\em Elec. J. Probab.} {\bf 19}, No. 16, 1-25, 2014.

\bibitem{TS2015}
T.~Simon. Mittag-Leffler functions and complete monotonicity. {\em Int. Transf. Spec. Funct.} {\bf 26} (1), 36-50, 2015.

\bibitem{T1939}
E.~C.~Titchmarsh. {\em The theory of functions.} Oxford University Press, 
Oxford, 1939.

\bibitem{Wiki}
Wikipedia, Mittag-Leffler function. {\tt https://en.wikipedia.org/wiki/Mittag-Leffler\_function}

\end{thebibliography}
\end{document}